\newif\ifarxiv\arxivfalse
\newif\ifaccel\accelfalse
\newif\ifcolors\colorsfalse
\arxivtrue

\ifarxiv
	\documentclass[10pt, twoside, reqno]{amsart}
\else
	\RequirePackage{amsmath}
	\documentclass[smallextended]{svjour3}
\fi

	\PassOptionsToPackage{noend}{algpseudocode}
	\PassOptionsToPackage{bookmarksdepth=4}{hyperref}
	
	\ifarxiv\else
		\makeatletter
			\let\cl@chapter\undefined
		\makeatother

		\smartqed
	\fi
	\usepackage{bm}
	\usepackage{mathrsfs}
	\usepackage[%
		eqreset=section,
		thmreset=section,
		noload={caption,subcaption},
	]{myPreamble}
	
	\makeatletter
		\def\operator@font{\rm}
	\makeatother
	\allowdisplaybreaks
	
	\ifarxiv\else
		
	\fi
	
	\ifcolors
		\def\red{\color{red}}
		\def\blue{\color{blue}}
	\else
		\def\red{}
		\def\blue{}
	\fi
	
	\renewcommand\FBsmooth{F}
	\renewcommand\FBnonsmooth{G}
	\renewcommand\FBstepsize{\Gamma}
	
	\AtBeginEnvironment{algorithmic}{%
		\setlength\baselineskip{1.3\baselineskip}%
	}
	
	\renewcommand\thealgorithm{\arabic{algorithm}}
	
	\newcommand\LL{{\mathscr L}_{\Gamma^{-1}}}
	\newcommand\ki[1][i]{t_\nu(#1)}
	
	\makeatletter
		\renewcommand{\appendixproof@toc}{subsubsection}
	
		\newcommand{\FBE}{\@ifstar\@@FBE\@FBE}
		\newcommand{\@FBE}{\Phi_\Gamma^{\text{\sc fb}}}
		\newcommand{\@@FBE}{\Phi_\gamma^{\text{\sc fb}}}
	
		\newcommand{\FBEC}{\tilde\Phi_\Gamma^{\text{\sc fb}}}
		
		\newcommand{\Res}{\@ifstar\@@Res\@Res}
		\newcommand{\@Res}{\operatorname R_\Gamma^{\text{\sc fb}}}
		\newcommand{\@@Res}{\operatorname R_\gamma^{\text{\sc fb}}}
	
		\newcommand{\T}{\@ifstar\@@T\@T}
		\newcommand{\@T}{\operatorname T_\Gamma^{\text{\sc fb}}}
		\newcommand{\@@T}{\operatorname T_\gamma^{\text{\sc fb}}}
	
		\newcommand{\E}{\@ifstar\@@E\@E}
		\newcommandx{\@E}[3][1=k,3={}]{
			\mathbb E_{#1}\ifstrempty{#2}{}{
				\left[
					#2\ifstrempty{#3}{}{
						\vphantom{#3}\right|\left.#3
					}
				\right]
			}
		}
		\newcommandx{\@@E}[3][1=k,3={}]{
			\mathbb E_{#1}\ifstrempty{#2}{}{
				[#2\ifstrempty{#3}{}{|#3}]
			}
		}
		
		\let\P\relax
		\newcommand{\P}{\@ifstar\@@P\@P}
		\newcommand{\@P}[2][k]{\mathcal P_{#1}\ifstrempty{#2}{}{\left[#2\right]}}
		\newcommand{\@@P}[2][k]{\mathcal P_{#1}\ifstrempty{#2}{}{[#2]}}
	\makeatother

		\newcommand\M[1][\Gamma]{\mathcal M_{#1}}
		\DeclareMathOperator{\blockdiag}{blkdiag}
	
		\let\mod\relax
		\DeclareMathOperator\mod{mod}
		\newcommand{\newnorm}[1]{
			\left\vert\kern-0.15ex\left\vert\kern-0.15ex\left\vert #1
			\right\vert\kern-0.15ex\right\vert\kern-0.15ex\right\vert
		}

\newcommand{\TheShortTitle}{Block-coordinate and incremental aggregated proximal gradient methods}
\newcommand{\TheTitle}{Block-coordinate and incremental aggregated proximal gradient methods for nonsmooth nonconvex problems}
\newcommand{\TheShortAuthor}{P. Latafat, A. Themelis and P. Patrinos}
\newcommand{\TheFunding}{%
	This work was supported by the Research Foundation Flanders (FWO) PhD grant 1196818N and research projects G086518N and G086318N;
	Research Council KU Leuven C1 project No. C14/18/068;
	Fonds de la Recherche Scientifique -- FNRS and the Fonds Wetenschappelijk Onderzoek -- Vlaanderen under EOS project no 30468160 (SeLMA)%
}
\newcommand{\TheKeywords}{%
	Nonsmooth nonconvex optimization,
	block-coordinate updates,
	forward-backward envelope,
	KL inequality%
}
\newcommand{\TheSubjclass}{%
	90C06, 
	90C25, 
	90C26, 
	49J52, 
	49J53.
}
\newcommand{\TheAbstract}{%
	This paper analyzes block-coordinate proximal gradient methods for minimizing the sum of a separable smooth function and a (nonseparable) nonsmooth function, both of which are allowed to be nonconvex.
	The main tool in our analysis is the forward-backward envelope (FBE), which serves as a particularly suitable continuous and real-valued Lyapunov function.
	Global and linear convergence results are established when the cost function satisfies the Kurdyka-\L ojasiewicz property without imposing convexity requirements on the smooth function.
	Two prominent special cases of the investigated setting are regularized finite sum minimization and the sharing problem; in particular, an immediate byproduct of our analysis leads to novel convergence results and rates for the popular Finito/MISO algorithm in the nonsmooth and nonconvex setting with very general sampling strategies.
	\ifaccel
		Moreover, favorable properties of the FBE under extra assumptions lead to an accelerated variant.
	\fi
}

\ifarxiv
	\title[\TheShortTitle]{\Large\scshape\bfseries\TheTitle}
	\author[\TheShortAuthor]{%
		Puya Latafat,
		Andreas Themelis and 
		Panagiotis Patrinos%
	}
	\thanks{%
		\TheAddressKU.
		{\tt
			\{%
				\href{mailto:puya.latafat@esat.kuleuven.be}{puya.latafat},%
				\href{mailto:andreas.themelis@esat.kuleuven.be}{andreas.themelis},%
				\href{mailto:panos.patrinos@esat.kuleuven.be}{panos.patrinos}%
			\}%
			\href{mailto:puya.latafat@esat.kuleuven.be,andreas.themelis@esat.kuleuven.be,panos.patrinos@esat.kuleuven.be}{@esat.kuleuven.be}%
		}%
	\\
		\TheFunding
	}
	\subjclass{\TheSubjclass}

\begin{document}

		\begin{abstract}
			\TheAbstract
		\end{abstract}

		\maketitle

\else

	\begin{document}
		\journalname{Mathematical Programming}

		\title{\TheTitle\thanks{\TheFunding.}}
		\titlerunning{\TheShortTitle}

		\author{%
			Puya Latafat\and
			Andreas Themelis\and
			Panagiotis Patrinos%
		}
		\authorrunning{\TheShortAuthor}

		\institute{%
			P. Latafat
			\at
			Tel.: +32 (0)16 374408\\
			\email{puya.latafat@kuleuven.be}
		\and
			A. Themelis
			\at
			Tel.: +32 (0)16 374573\\
			\email{andreas.themelis@kuleuven.be}
		\and
			P. Patrinos
		\at
			Tel.: +32 (0)16 374445\\
			\email{panos.patrinos@esat.kuleuven.be}%
		\and
			\TheAddressKU.
		}

		\date{Received: date / Accepted: date}

		\maketitle

		\begin{abstract}
			\TheAbstract
			\subclass{\TheSubjclass}%
		\end{abstract}
\fi


	\section{Introduction}\label{sec:Introduction}

	This paper addresses block-coordinate (BC) proximal gradient methods for problems of the form
	\begin{equation}\label{eq:P}
		\minimize_{\bm x=(x_1,\dots,x_N)\in\R^{\sum_in_i}}
			\Phi(\bm x)
		{}\coloneqq{}
			F(\bm x)
			{}+{}
			G(\bm x),
		\quad\text{where}\quad
		\textstyle
		F(\bm x)\coloneqq\tfrac1N\sum_{i=1}^N f_i(x_i),
	\end{equation}
	in the following setting.
	
	\begin{ass}[problem setting]\label{ass:basic}%
		In problem \eqref{eq:P} the following hold:
		\begin{enumeratass}
		\item\label{ass:f}%
			function \(f_i\) is \(L_{f_i}\)-smooth (Lipschitz differentiable with modulus \(L_{f_i}\)), \(i\in[N]\);
		\item\label{ass:g}%
			function \(G\) is proper and lower semicontinuous (lsc);
		\item\label{ass:phi}%
			a solution exists: \(\argmin\Phi\neq\emptyset\).
		\end{enumeratass}
	\end{ass}
	
	Unlike typical cases analyzed in the literature where \(G\) is separable \cite{tseng2001convergence,tseng2009coordinate,nesterov2012efficiency,beck2013convergence,bolte2014proximal,richtarik2014iteration,lin2015accelerated,chouzenoux2016block,hong2017iteration,xu2017globally}, we here consider the complementary case where it is only the smooth term \(F\) that is assumed to be separable. The main challenge in analyzing convergence of BC schemes for \eqref{eq:P} especially in the nonconvex setting is the fact that even in expectation the cost does not necessarily decrease along the trajectories.
	Instead, we demonstrate that the  forward-backward envelope (FBE) \cite{patrinos2013proximal,themelis2018forward} is a suitable Lyapunov function for such problems.  
	
	Several BC-type algorithms that allow for a nonseparable nonsmooth term have been considered in the literature, however, all in convex settings.  
	In \cite{tseng2008block,tseng2010coordinate} a class of convex composite problems is studied that involves a linear constraint as the nonsmooth nonseparable term. 
	A BC algorithm with a Gauss-Southwell-type rule is proposed and the convergence is established using the cost as Lyapunov function by exploiting linearity of the constraint to ensure feasibility. A refined analysis in \cite{necoara2013random,necoara2014random} extends this to a random coordinate selection strategy. 
	Another approach in the convex case is to consider randomized BC updates applied to general averaged operators. Although this approach can allow for fully nonseparable problems, usually  separable nonsmooth functions are considered in the literature. 
	The convergence analysis of such methods relies on establishing quasi-Fej\'er monotonicity \cite{iutzeler2013asynchronous,combettes2015stochastic,pesquet2015class,bianchi2016coordinate,peng2016arock,latafat2019new}.
	In a primal-dual setting in \cite{fercoq2019coordinate} a combination of Bregman and Euclidean distance is employed as Lyapunov function.
	In \cite{hanzely2018sega} a BC algorithm is proposed for strongly convex algorithms that involves  coordinate updates for the gradient followed by a full proximal step, and the distance from the (unique) solution is used as Lyapunov function.
	The analysis and the Lyapunov functions in all of the above mentioned works rely heavily on convexity and are not suitable for nonconvex settings.  
	
	Thanks to the nonconvexity and nonseparability of \(G\), many machine learning problems can be formulated as in \eqref{eq:P}, a primary example being constrained and/or regularized finite sum problems \cite{bertsekas2011incremental,shalevshwartz2013stochastic,defazio2014finito,defazio2014saga,mairal2015incremental,reddi2016proximal,reddi2016stochastic,schmidt2017minimizing}
	\begin{equation}\label{eq:FSP}
	\textstyle
		\minimize_{x\in\R^n}
			\varphi(x)
		{}\coloneqq{}
			\tfrac1N\sum_{i=1}^N f_i(x)
			{}+{}
			g(x),
	\end{equation}
	where \(\func{f_i}{\R^n}{\R}\) are smooth functions and \(\func{g}{\R^n}{\Rinf}\) is possibly nonsmooth, and everything here can be nonconvex.
	In fact, one way to cast \eqref{eq:FSP} into the form of problem \eqref{eq:P} is by setting
	\begin{equation}\label{eq:FINITOG}
	\textstyle
		G(\bm x)
	{}\coloneqq{}
		\tfrac1N\sum_{i=1}^Ng(x_i)
		{}+{}
		\indicator_C(\bm x),
	\end{equation}
	where
	\(
		C
	{}\coloneqq{}
		\set{\bm x\in\R^{nN}}[x_1=x_2=\dots=x_N]
	\)
	is the consensus set, and \(\indicator_C\) is the indicator function of set \(C\), namely
	\(
		\indicator_C(\bm x)=0
	\)
	for \(\bm x\in C\) and \(\infty\) otherwise.
	Since the nonsmooth term \(g\) is allowed to be nonconvex, formulation \eqref{eq:FSP} can account for nonconvex constraints such as rank constraints or zero norm balls, and nonconvex regularizers such as \(\ell^p\) with \(p\in[0,1)\), \cite{hou2012complexity}. 
	
	Another prominent example in distributed applications is the \emph{``sharing''} problem \cite{boyd2011distributed}:
	\begin{equation}\label{eq:SP}
		\minimize_{\bm x\in\R^{nN}}\Phi(\bm x)
	{}\coloneqq{}
		\textstyle
		\tfrac1N\sum_{i=1}^Nf_i(x_i)
		{}+{}
		g\Bigl(\sum_{i=1}^Nx_i\Bigr)
		.
	\end{equation}
	where  \(\func{f_i}{\R^n}{\R}\) are smooth functions and \(\func{g}{\R^n}{\Rinf}\) is nonsmooth, and all are possibly nonconvex. The sharing problem is cast as in \eqref{eq:P} by setting \(G\coloneqq g \circ A\), where \(A\coloneqq[\I_n~\dots~\I_n]\in\R^{n\times nN}\) (\(I_r\) denotes the \(r\times r\) identity matrix).

	\subsection{The main block-coordinate algorithm}\label{sec:BC}

	While gradient evaluations are the building blocks of smooth minimization, a fundamental tool to deal with a nonsmooth lsc term \(\func{\psi}{\R^r}{\Rinf}\) is its \DEF{\(V\)-proximal mapping}
	\begin{equation}\label{eq:prox}
		\prox_\psi^V(x)
	{}\coloneqq{}
		\argmin_{w\in\R^r}\set{
			\psi(w)
			{}+{}
			\tfrac12\|w-x\|^2_V
		},
	\end{equation}
	where \(V\) is a symmetric and positive definite matrix and \(\|{}\cdot{}\|_V\) indicates the norm induced by the scalar product \((x,y)\mapsto\innprod{x}{Vy}\).
	It is common to take \(V=t^{-1}\I_r\) as a multiple of the \(r\times r\) identity matrix \(\I_r\), in which case the notation \(\prox_{t\psi}\) is typically used and \(t\) is referred to as a stepsize.
	While this operator enjoys nice regularity properties when \(g\) is convex, such as (single valuedness and) Lipschitz continuity, for nonconvex \(g\) it may fail to be a well-defined function and rather has to be intended as a point-to-set mapping \(\ffunc{\prox_\psi^V}{\R^r}{\R^r}\).
	Nevertheless, the value function associated to the minimization problem in the definition \eqref{eq:prox}, namely the \emph{Moreau envelope}
	\begin{equation}\label{eq:Moreau}
		\psi^V(x)
	{}\coloneqq{}
		\min_{w\in\R^r}\set{
			\psi(w)
			{}+{}
			\tfrac12\|w-x\|^2_V
		},
	\end{equation}
	is a well-defined real-valued function, in fact locally Lipschitz continuous, that lower bounds \(\psi\) and shares with \(\psi\) infima and minimizers.
	The proximal mapping is available in closed form for many useful functions, many of which are widely used regularizers in machine learning; for instance, the proximal mapping of the \(\ell^0\) and \(\ell^1\) regularizers amount to hard and soft thresholding operators.
	
	In many applications the cost to be minimized is structured as the sum of a smooth term \(h\) and a proximable (\ie with easily computable proximal mapping) term \(\psi\).
	In these cases, the \emph{proximal gradient method} \cite{fukushima1981generalized,attouch2013convergence} constitutes a cornerstone iterative method that interleaves gradient descent steps on the smooth function and proximal operations on the nonsmooth function, resulting in iterations of the form
	\(
		x^+
	{}\in{}
		\prox_{\gamma\psi}(x-\gamma\nabla h(x))
	\)
	for some suitable stepsize $\gamma$.
	
	Our proposed scheme to address problem \eqref{eq:P} is a BC variant of the proximal gradient method, in the sense that only some coordinates are updated according to the proximal gradient rule, while the others are left unchanged.
	This concept is synopsized in \Cref{alg:BC}, which constitutes the general algorithm addressed in this paper.
	
	\begin{algorithm}
		\caption{General forward-backward block-coordinate scheme}
		\label{alg:BC}
	\begin{algorithmic}[1]
	\Require
		\(\bm x^0\in\R^{\sum_in_i}\),~
		\(\gamma_i\in(0,\nicefrac{N}{L_{f_i}})\),
		{\small \(i\in[N]\)}
	\Statex 
		\(
			\Gamma=\blockdiag(\gamma_1\I_{n_1},\dots,\gamma_N\I_{n_N})
		\),~
		\(k=0\)
	\item[{\sc Repeat} until convergence]
		\State
			\(
				\bm z^k
			{}\in{}
				\prox_G^{\Gamma^{-1}}\bigl(
					\bm x^k-\Gamma\nabla F(\bm x^k)
				\bigr)
			\)
		\State\label{state:BC:sampling}%
			select a set of indices \(I^{k+1}\subseteq[N]\)
		\State
			update~~
			\(x_i^{k+1}= z_i^k\)
			~for \(i\in I^{k+1}\)
			~~and~~
			\(x_i^{k+1}= x_i^k\)
			~for \(i\notin I^{k+1}\),~
			\(k\gets k+1\)
	\item[{\sc Return}]
		\(\bm z^k\)
	\end{algorithmic}
	\end{algorithm}
	
	Although seemingly wasteful, in many cases one can efficiently compute individual blocks without the need of full operations.
	In fact BC \Cref{alg:BC} bridges the gap between a BC framework and a class of incremental methods where a global computation typically involving the full gradient is carried out incrementally via performing computations only for a subset of coordinates. 
	Two such broad applications, problems \eqref{eq:FSP} and \eqref{eq:SP}, are discussed in the dedicated \Cref{sec:Finito,sec:Sharing}, where among other things we will show that \Cref{alg:BC} leads to the well known Finito/MISO algorithm \cite{defazio2014finito,mairal2015incremental}.

		\subsection{Contribution}

	\begin{enumerate}[%
		leftmargin=0pt,
		labelwidth=7pt,
		itemindent=\labelwidth+\labelsep,
		label=\rlap{{\bf\arabic*)}}\hspace*{\labelwidth},
	]
	\item 
		To the best of our knowledge this is the first analysis of BC schemes with a nonseparable nonsmooth term and in the fully nonconvex setting.
		While the original cost \(\Phi\) cannot serve as a Lyapunov function, we show that the forward-backward envelope (FBE) \cite{patrinos2013proximal,themelis2018forward} decreases surely, not only in expectation (\Cref{thm:sure}).
	\item
		This allows for a quite general convergence analysis for different sampling criteria. 
		This paper in particular covers randomized strategies (\Cref{sec:random}) where at each iteration one or more coordinates are sampled with possibly time-varying probabilities, as well as essentially cyclic (and in particular cyclic and shuffled) strategies in case the nonsmooth term is convex (\Cref{sec:cyclic}).
	\item
		We exploit the Kurdyka-\L ojasiewicz (KL) property to show global (as opposed to subsequential) and linear convergence when the sampling is essentially cyclic and the nonsmooth function is convex, without imposing convexity requirements on the smooth functions (\Cref{thm:cyclic:global}).
	\ifaccel
		\item
			When \(G\) is convex and \(F\) is twice continuously differentiable, the FBE is continuously differentiable.
			If, additionally, \(F\) is (strongly) convex and quadratic, then the FBE is (strongly) convex and has Lipschitz-continuous gradient.
			Owing to these favorable properties, we propose a new BC Nesterov-type acceleration algorithm for minimizing the sum of a block-separable convex quadratic plus a nonsmooth convex function, whose analysis directly follows from existing work on smooth BC minimization \cite{allen2016even}.
	\fi
	\item
		As immediate byproducts of our analysis we obtain
		{\bf (a)} an incremental algorithm for the sharing problem \cite{boyd2011distributed} that to the best of our knowledge is novel (\Cref{sec:Sharing}), and
		{\bf (b)} the Finito/MISO algorithm \cite{defazio2014finito,mairal2015incremental} leading to a much simpler and more general analysis than available in the literature with new convergence results both for randomized sampling strategies in the fully nonconvex setting and for essentially cyclic samplings when the nonsmooth term is convex (\Cref{sec:Finito}).
	\end{enumerate}

		\subsection{Organization}

	The rest of the paper is organized as follows.
	The core of the paper lies in the convergence analysis of \Cref{alg:BC} detailed in \Cref{sec:convergence}: \Cref{sec:FBE} introduces the FBE, fundamental tool of our methodology and lists some of its properties whose proofs are detailed in the dedicated \Cref{sec:proofs:FBE}, followed by other ancillary results documented in \Cref{sec:auxiliary}.
	The algorithmic analysis begins in \Cref{sec:sure} with a collection of facts that hold independently of the chosen sampling strategy, and later specializes to randomized and essentially cyclic samplings in the dedicated \Cref{sec:random,sec:cyclic}.
	\Cref{sec:Finito,sec:Sharing} discuss two particular instances of the investigated algorithmic framework, namely (a generalization of) the Finito/MISO algorithm for finite sum minimization and an incremental scheme for the sharing problem, both for fully nonconvex and nonsmooth formulations.
	Convergence results are immediately inferred from those of the more general BC \Cref{alg:BC}.
	\Cref{sec:Conclusions} concludes the paper.

	\section{Convergence analysis}\label{sec:convergence}

	We begin by observing that \Cref{ass:basic} is enough to guarantee the well definedness of the forward-backward operator in \Cref{alg:BC}, which for notational convenience will be henceforth denoted as \(\T(\bm x)\).
	Namely, \(\ffunc{\T}{\R^{\sum_in_i}}{\R^{\sum_in_i}}\) is the point-to-set mapping
	\begin{align*}
		\T(\bm x)
	{}\coloneqq{} &
		\prox_G^{\Gamma^{-1}}\left(\Fw{\bm x}\right)
	\\
	\numberthis\label{eq:T}
	{}={} &
		\argmin_{\bm w\in\R^{\sum_in_i}}\set{
			F(\bm x)+\innprod{\nabla F(\bm x)}{\bm w-\bm x}
			{}+{}
			G(\bm w)
			{}+{}
			\tfrac12\|\bm w-\bm x\|_{\Gamma^{-1}}^2
		}.
	\end{align*}
	
	\begin{lem}\label{thm:osc}%
		Suppose that \Cref{ass:basic} holds, and let \(\Gamma\coloneqq\blockdiag(\gamma_1\I_{n_1},\dots,\gamma_N\I_{n_N})\) with \(\gamma_i\in(0,\nicefrac{N}{L_{f_i}})\), \(i\in[N]\).
		Then \(\prox_G^{\Gamma^{-1}}\) and \(\T\) are locally bounded, outer semicontinuous (osc), nonempty- and compact-valued mappings.
		\begin{proof}
			See \Cref{proof:thm:osc}.
		\end{proof}
	\end{lem}

		\subsection{The forward-backward envelope}\label{sec:FBE}

	The fundamental challenge in the analysis of \eqref{eq:P} is the fact that, without separability of \(G\), descent on the cost function cannot be established even in expectation.
	Instead, we show that the \emph{forward-backward envelope} (FBE) \cite{patrinos2013proximal,themelis2018forward} can be used as Lyapunov function.
	This subsection formally introduces the FBE, here generalized to account for a matrix-valued stepsize parameter \(\Gamma\), and lists some of its basic properties needed for the convergence analysis of \Cref{alg:BC}.
	Although easy adaptations of the similar results in \cite{patrinos2013proximal,themelis2018forward,themelis2019acceleration}, for the sake of self-inclusiveness the proofs are detailed in the dedicated \Cref{sec:proofs:FBE}.
	
	\begin{subequations}
		\begin{defin}[forward-backward envelope]\label{def:FBE}%
			In problem \eqref{eq:P}, let \(f_i\) be differentiable functions, \(i\in[N]\), and for \(\gamma_1,\dots,\gamma_N>0\) let
			\(
				\Gamma=\blockdiag(\gamma_1\I_{n_1},\dots,\gamma_N\I_{n_N})
			\).
			The forward-backward envelope (FBE) associated to \eqref{eq:P} with stepsize \(\Gamma\) is the function
			\(
				\func{\FBE}{\R^{\sum_in_i}}{[-\infty,\infty)}
			\)
			defined as
			\begin{equation}
				\label{eq:FBE}
				\FBE(\bm x)
			{}\coloneqq{}
				\inf_{\bm w\in\R^{\sum_in_i}}\set{
					F(\bm x)+\innprod{\nabla F(\bm x)}{\bm w-\bm x}
					{}+{}
					G(\bm w)
					{}+{}
					\tfrac12\|\bm w-\bm x\|_{\Gamma^{-1}}^2
				}.
			\end{equation}
		\end{defin}
	
		\Cref{def:FBE} highlights an important symmetry between the Moreau envelope and the FBE: similarly to the relation between the Moreau envelope \eqref{eq:Moreau} and the proximal mapping \eqref{eq:prox}, the FBE \eqref{eq:FBE} is the value function associated with the proximal gradient mapping \eqref{eq:T}.
		By replacing any minimizer \(\bm z\in\T(\bm x)\) in the right-hand side of \eqref{eq:FBE} one obtains yet another interesting interpretation of the FBE in terms of the \(\Gamma^{-1}\)-augmented Lagrangian associated to \eqref{eq:P}
		\begin{align}
			\nonumber
			\LL(\bm x,\bm z,\bm y)
		{}\coloneqq{} &
			F(\bm x)+G(\bm z)+\innprod{\bm y}{\bm x-\bm z}
			{}+{}
			\tfrac12\|\bm x-\bm z\|_{\Gamma^{-1}}^2,
		\shortintertext{namely,}
			\label{eq:FBEz}
			\FBE(\bm x)
		{}={} &
			F(\bm x)+\innprod{\nabla F(\bm x)}{\bm z-\bm x}
			{}+{}
			G(\bm z)
			{}+{}
			\tfrac12\|\bm z-\bm x\|_{\Gamma^{-1}}^2
		\\
		{}={} &
			\LL(\bm x,\bm z,-\nabla F(\bm x)).
		\shortintertext{%
			Lastly, by rearranging the terms it can easily be seen that
		}
		\label{eq:FBEMoreau}
			\FBE(\bm x)
		{}={} &
			F(\bm x)
			{}-{}
			\tfrac12\|\nabla F(\bm x)\|_\Gamma^2
			{}+{}
			G^{\Gamma^{-1}}(\Fw{\bm x}),
		\end{align}
		hence in particular that the FBE inherits regularity properties of \(G^{\Gamma^{-1}}\) and \(\nabla F\), some of which are summarized in the next result.
	\end{subequations}
	
	\begin{lem}[FBE: fundamental inequalities]\label{thm:FBEineq}%
		Suppose that \Cref{ass:basic} is satisfied and let \(\gamma_i\in(0,\nicefrac{N}{L_{f_i}})\), \(i\in[N]\).
		Then, the FBE \(\FBE\) is a (real-valued and) locally Lipschitz-continuous function.
		Moreover, the following hold for any \(\bm x\in\R^{\sum_in_i}\):
		\begin{enumerate}
		\item\label{thm:leq}
			\(\FBE(\bm x)\leq\Phi(\bm x)\).
		\item\label{thm:geq}
			\(
				\tfrac12\|\bm z-\bm x\|^2_{\Gamma^{-1}-\Lambda_F}
			{}\leq{}
				\FBE(\bm x)-\Phi(\bm z)
			{}\leq{}
				\tfrac12\|\bm z-\bm x\|^2_{\Gamma^{-1}+\Lambda_F}
			\)
			for any \(\bm z\in\T(\bm x)\), where
			\(
				\Lambda_F
			{}\coloneqq{}
				\tfrac1N
				\blockdiag\bigl(L_{f_1}\I_{n_1},\dots, L_{f_n}\I_{n_N}\bigr)
			\).
		\item\label{thm:strconcost}
			If in addition each $f_i$ is $\mu_{f_i}$-strongly convex and $G$ is convex, then for every \(\bm x\in\R^{\sum_in_i}\)
			\[
				\tfrac12\|\bm z-\bm x^\star\|_{\mu_F}^2
			{}\leq{}
				\FBE(\bm x)-\min\Phi
			\]
			where \(\bm x^\star\coloneqq\argmin\Phi\),
			\(
				\mu_F
			{}\coloneqq{}
				\frac1N\blockdiag\bigl(\mu_{f_1}\I_{n_1},\dots,\mu_{f_N}\I_{n_N}\bigr)
			\),
			and \(\bm z=\T(\bm x)\).
		\end{enumerate}
		\begin{proof}
			See \Cref{proof:thm:FBEineq}.
		\end{proof}
	\end{lem}
	
	Another key property that the FBE shares with the Moreau envelope is that minimizing the extended-real valued function \(\Phi\) is equivalent to minimizing the continuous function \(\FBE\).
	Moreover, the former is level bounded iff so is the latter.
	This fact will be particularly useful for the analysis of \Cref{alg:BC}, as it will be shown in \Cref{thm:sure} that the FBE (surely) decreases along its iterates.
	As a consequence, despite the fact that the same does not hold for \(\Phi\) (in fact, iterates may even be infeasible), coercivity of \(\Phi\) is enough to guarantee boundedness of \(\seq{\bm x^k}\) and \(\seq{\bm z^k}\).
	
	\begin{lem}[FBE: minimization equivalence]\label{thm:FBEmin}%
		Suppose that \Cref{ass:basic} is satisfied and that \(\gamma_i\in(0,\nicefrac{N}{L_i})\), \(i\in[N]\).
		Then the following hold:
		\begin{enumerate}
		\item\label{thm:min}
			\(\min\FBE=\min\Phi\);
		\item\label{thm:argmin}
			\(\argmin\FBE=\argmin\Phi\);
		\item\label{thm:LB}
			\(\FBE\) is level bounded iff so is \(\Phi\).
		\end{enumerate}
		\begin{proof}
			See \Cref{proof:thm:FBEmin}.
		\end{proof}
	\end{lem}
	
	We remark that the kinship of \(\FBE\) and \(\Phi\) extends also to local minimality; the interested reader is referred to \cite[Th. 3.6]{themelis2018proximal} for details.

		\subsection{A sure descent lemma}\label{sec:sure}

	We now proceed to the theoretical analysis of \Cref{alg:BC}.
	Clearly, some assumptions on the index selection criterion are needed in order to establish reasonable convergence results, for little can be guaranteed if, for instance, one of the indices is never selected.
	Nevertheless, for the sake of a general analysis it is instrumental to first investigate which properties hold independently of such criteria.
	After listing some of these facts in \Cref{thm:sure}, in \Cref{sec:random,sec:cyclic} we will specialize the results to randomized and (essentially) cyclic sampling strategies.
	
	\begin{lem}[sure descent]\label{thm:sure}%
		Suppose that \Cref{ass:basic} is satisfied.
		Then, the following hold for the iterates generated by \Cref{alg:BC}:
		\begin{enumerate}
		\item\label{thm:Igeq}
			\(
				\FBE(\bm x^{k+1})
			{}\leq{}
				\FBE(\bm x^k)
				{}-{}
				\sum_{i\in I^{k+1}}\tfrac{\xi_i}{2\gamma_i}\|z_i^k-x_i^k\|^2
			\),
			where \(\xi_i\coloneqq\frac{N-\gamma_iL_{f_i}}{N}\), \(i\in[N]\), are strictly positive;
		\item\label{thm:decrease}%
			\(\seq{\FBE(\bm x^k)}\) monotonically decreases to a finite value \(\Phi_\star\geq\min\Phi\);
		\item\label{thm:omega}%
			\(\FBE\) is constant (and equals \(\Phi_\star\) as above) on the set of accumulation points of \(\seq{\bm x^k}\);
		\item\label{thm:xdiff}%
			the sequence \(\seq{\|\bm x^{k+1}-\bm x^k\|^2}\) has finite sum (and in particular vanishes);
		\item\label{thm:bounded}%
			if \(\Phi\) is coercive, then \(\seq{\bm x^k}\) and \(\seq{\bm z^k}\) are bounded.
		\end{enumerate}
		\begin{proof}
	\begin{proofitemize}
	\item\ref{thm:Igeq}~
		To ease notation, let
		\(
			\Lambda_F
		{}\coloneqq{} \tfrac1N
			\blockdiag\bigl(L_{f_1}\I_{n_1},\dots, L_{f_n}\I_{n_N}\bigr)
		\) 
		and for \(\bm w\in\R^{\sum_in_i}\) let \(w_I\in\R^{\sum_{i\in I}n_i}\) denote the slice \((w_i)_{i\in I}\), and let \(\Lambda_{F_I},\Gamma_I\in\R^{\sum_{i\in I}n_i\times\sum_{i\in I}n_i}\) be defined accordingly.
		Start by observing that, since \(\bm z^{k+1}\in\prox_G^{\Gamma^{-1}}(\Fw{\bm x^{k+1}})\), from the proximal inequality on $G$ it follows that
		\begin{align*}
			G(\bm z^{k+1})-G(\bm z^k)
		{}\leq{} &
			\tfrac12\|\bm z^k-\bm x^{k+1}+\Gamma\nabla F(\bm x^{k+1})\|_{\Gamma^{-1}}^2
			{}-{}
			\tfrac12\|\bm z^{k+1}-\bm x^{k+1}+\Gamma\nabla F(\bm x^{k+1})\|_{\Gamma^{-1}}^2
		\\
		={} &
		\numberthis\label{eq:proxIneqFBS}
			\tfrac12\|\bm z^k-\bm x^{k+1}\|_{\Gamma^{-1}}^2
			{}-{}
			\tfrac12\|\bm z^{k+1}-\bm x^{k+1}\|_{\Gamma^{-1}}^2
			{}+{}
			\innprod{\nabla F(\bm x^{k+1})}{\bm z^k-\bm z^{k+1}}.
		\end{align*}
		We have
		\ifarxiv\else
		\bgroup\mathtight[0.5]%
		\fi
		\begin{align*}
			\FBE(\bm x^{k+1})-\FBE(\bm x^k)
		{}={} &
			{\red
				F(\bm x^{k+1})
			}
			{}+{}
			\innprod{\nabla F(\bm x^{k+1})}{\bm z^{k+1}-\bm x^{k+1}}
			{\blue
				{}+{}
				G(\bm z^{k+1})
			}
			{}+{}
			\tfrac12\|\bm z^{k+1}-\bm x^{k+1}\|_{\Gamma^{-1}}^2
		\\
		&
			{}-{}
			\left(
				{\red
					F(\bm x^k)+\innprod{\nabla F(\bm x^k)}{\bm z^k-\bm x^k}
				}
				{\blue
					{}+{}
					G(\bm z^k)
				}
				{}+{}
				\tfrac12\|\bm z^k-\bm x^k\|_{\Gamma^{-1}}^2
			\right)
		\shortintertext{%
			{\red
				apply the upper bound in \eqref{eq:Lip} with \(\bm w=\bm x^{k+1}\)
			}%
			and
			{\blue
				the proximal inequality \eqref{eq:proxIneqFBS}
			}%
		}
		{}\leq{} &
			{\red
				\innprod{\nabla F(\bm x^k)}{\bm x^{k+1}-\bm z^k}
				{}+{}
				\tfrac12\|\bm x^{k+1}-\bm x^k\|_{\Lambda_F}^2
			}
			{}+{}
			\innprod{\nabla F(\bm x^{k+1})}{{\blue\bm z^k}-\bm x^{k+1}}
		\\
		&
			{}-{}
			\tfrac12\|\bm z^k-\bm x^k\|_{\Gamma^{-1}}^2
			{\blue
				{}+{}
				\tfrac12\|\bm z^k-\bm x^{k+1}\|_{\Gamma^{-1}}^2
			}.
		\end{align*}
		\ifarxiv\else
			\egroup
		\fi
		To conclude, notice that the \(\ell\)-th block of \(\nabla F(\bm x^k)-\nabla F(\bm x^{k+1})\) is zero for \(\ell\notin I\), and that the \(\ell\)-th block of \(\bm x^{k+1}-\bm z^k\) is zero if \(\ell\in I\).
		Hence, the scalar product vanishes.
		For similar reasons, one has
		\(
			\|
				\bm z^k-\bm x^{k+1}
			\|^2_{\Gamma^{-1}}
			{}-{}
			\|\bm z^k-\bm x^k\|_{\Gamma^{-1}}^2
		{}={}
			{}-{}
			\|z_I^k-x_I^k\|_{\Gamma_I^{-1}}^2
		\)
		and
		\(
			\|\bm x^{k+1}-\bm x^k\|_{\Lambda_F}^2
		{}={}
			\|z_I^k-x_I^k\|_{\Lambda_{F_I}}^2
		\),
		yielding the claimed expression.
	\item\ref{thm:decrease}~
		Monotonic decrease of \(\seq{\FBE(\bm x^k)}\) is a direct consequence of assert \ref{thm:Igeq}.
		This ensures that the sequence converges to some value \(\Phi_\star\), bounded below by \(\min\Phi\) in light of \Cref{thm:min}.
	\item\ref{thm:omega}~
		Directly follows from assert \ref{thm:decrease} together with the continuity of \(\FBE\), see \Cref{thm:FBEineq}.
	\item\ref{thm:xdiff}~
		Denoting
		\(
			\xi_{\rm min}
		{}\coloneqq{}
			\min_{i\in[N]}\set{
				\xi_i
			}
		\)
		which is a strictly positive constant, it follows from assert \ref{thm:Igeq} that for each \(k\in\N\) it holds that
		\begin{align*}
			\FBE(\bm x^{k+1})-\FBE(\bm x^k)
		{}\leq{} &
			{}-{}
			\sum_{\mathclap{i\in I^{k+1}}}{
				\tfrac{\xi_i}{2\gamma_i}\|z_i^k-x_i^k\|^2
			}
		\\
		{}\leq{} &
			{}-{}
			\tfrac{\xi_{\rm min}}{2}
			\sum_{i\in I^{k+1}}{
				\gamma_i^{-1}\|z_i^k-x_i^k\|^2
			}
		\\
		\numberthis\label{eq:SDx}
		{}={} &
			{}-{}
			\tfrac{\xi_{\rm min}}{2}
			\|\bm x^{k+1}-\bm x^k\|_{\Gamma^{-1}}^2.
		\end{align*}
		By summing for \(k\in\N\) and using the positive definiteness of \(\Gamma^{-1}\) together with the fact that \(\min\FBE=\min\Phi>\infty\) as ensured by \Cref{thm:min} and \Cref{ass:phi}, we obtain that
		\(
			\sum_{k\in\N}\|\bm x^{k+1}-\bm x^k\|^2
		{}<{}
			\infty
		\).
	\item\ref{thm:bounded}~
		It follows from assert \ref{thm:decrease} that the entire sequence \(\seq{\bm x^k}\) is contained in the sublevel set \(\set{\bm w}[\FBE(\bm w)\leq\FBE(\bm x^0)]\), which is bounded provided that \(\Phi\) is coercive as shown in \Cref{thm:LB}.
		In turn, boundedness of \(\seq{\bm z^k}\) then follows from local boundedness of \(\T\), cf. \Cref{thm:osc}.
	\qedhere
	\end{proofitemize}
		\end{proof}
	\end{lem}

		\subsection{Randomized sampling}\label{sec:random}

	In this section we provide convergence results for \Cref{alg:BC} where the index selection criterion complies with the following requirement.
	
	\begin{ass}[randomized sampling requirements]\label{ass:random}%
		There exist \(p_1,\dots,p_N>0\) such that, at any iteration and independently of the past, each $i\in[N]$ is sampled with probability at least $p_i$.
	\end{ass}
	
	Our notion of randomization is general enough to allow for time-varying probabilities and mini-batch selections.
	The role of parameters \(p_i\) in \Cref{ass:random} is to prevent that an index is sampled with arbitrarily small probability.
	In more rigorous terms,
	\(
		\P*{i\in I^{k+1}}
	{}\geq{}
		p_i
	\)
	shall hold for all \(i\in[N]\), where \(\P{}\) represents the probability conditional to the knowledge at iteration \(k\).
	Notice that we do not require the \(p_i\)'s to sum up to one, as multiple index selections are allowed, similar to the setting of \cite{bianchi2016coordinate,latafat2019new} in the convex case.
	
	Due to the possible nonconvexity of problem \eqref{eq:P}, unless additional assumptions are made not much can be said about convergence of the iterates to a unique point.
	Nevertheless, the following result shows that any accumulation point \(\bm x^\star\) of sequences $\seq{\bm x^k}$ and $\seq{\bm z^k}$ generated by \Cref{alg:BC} is a stationary point, in the sense that it satisfies the necessary condition for minimality
	\(
		0\in\hat\partial\Phi(\bm x^\star)
	\),
	where \(\hat\partial\) denotes the (regular) nonconvex subdifferential, see \cite[Th. 10.1]{rockafellar2011variational}.
	
	\begin{thm}[randomized sampling: subsequential convergence]\label{thm:random:subseq}%
		Suppose that \Cref{ass:basic,ass:random} are satisfied.
		Then, the following hold almost surely for the iterates generated by \Cref{alg:BC}:
		\begin{enumerate}
		\item\label{thm:res}%
			the sequence \(\seq{\|\bm x^k-\bm z^k\|^2}\) has finite sum (and in particular vanishes);
		\item\label{thm:decreasez}%
			the sequence \(\seq{\Phi(\bm z^k)}\) converges to \(\Phi_\star\) as in \Cref{thm:decrease};
		\item\label{thm:cluster}%
			\(\seq{\bm x^k}\) and \(\seq{\bm z^k}\) have same cluster points, all stationary and on which \(\Phi\) and \(\FBE\) equal \(\Phi_\star\).%
		\end{enumerate}
		\begin{proof}
	
	In what follows, \(\E{}\) denotes the expectation conditional to the knowledge at iteration \(k\).
	\begin{proofitemize}
	\item\ref{thm:res}~
		Let
		\(
			\xi_i\coloneqq\frac{N-\gamma_iL_{f_i}}{N}>0
		\),
		\(i\in[N]\), be as in \Cref{thm:Igeq}.
		We have
		\begin{align*}
			\E{\FBE(\bm x^{k+1})}
		{}\overrel*[\leq]{\ref{thm:Igeq}}{} &
			\E{
				\FBE(\bm x^k)
				{}-{}
				\sum_{i\in I^{k+1}}{
					\tfrac{\xi_i}{2\gamma_i}\|z_i^k-x_i^k\|^2
				}
			}
		\\
		{}={} &
			\FBE(\bm x^k)
			{}-{}
			\sum_{I\in\Omega}{
				\P{\mathcal I^{k+1}=I}
				\sum_{i\in I}{
					\tfrac{\xi_i}{2\gamma_i}\|z_i^k-x_i^k\|^2
				}
			}
		\\
		{}={} &
			\FBE(\bm x^k)
			{}-{}
			\sum_{i=1}^N{
				\sum_{I\in\Omega,I\ni i}{
					\P{\mathcal I^{k+1}=I}
					\tfrac{\xi_i}{2\gamma_i}\|z_i^k-x_i^k\|^2
				}
			}
		\\
		\numberthis\label{eq:EFBE+}
		{}\leq{} &
			\FBE(\bm x^k)
			{}-{}
			\sum_{i=1}^N{
				\tfrac{p_i\xi_i}{2\gamma_i}\|z_i^k-x_i^k\|^2
			},
		\end{align*}
		where $\Omega\subseteq 2^{[N]}$ is the sample space ($2^{[N]}$ denotes the power set of $[N]$).  
		Therefore,
		\begin{equation}\label{eq:ExSD}
			\E{\FBE(\bm x^{k+1})}
		{}\leq{}
			\FBE(\bm x^k)
			{}-{}
			\tfrac\sigma2
			\|\bm x^k-\bm z^k\|_{\Gamma^{-1}}^2
		\quad\text{where }
			\sigma
		{}\coloneqq{}
			\min_{i=1\dots N}{
				p_i\xi_i
			}
			{}>{}
			0.
		\end{equation}
		The claim follows from the Robbins-Siegmund supermartingale theorem, see \eg \cite{robbins1985convergence} or \cite[Prop. 2]{bertsekas2011incremental}.
	\item\ref{thm:decreasez}~
		Observe that
		\(
			\FBE(\bm x^k)-\|\bm z^k-\bm x^k\|^2_{\Gamma^{-1}+\Lambda_F}
		{}\leq{}
			\Phi(\bm z^k)
		{}\leq{}
			\FBE(\bm x^k)-\|\bm z^k-\bm x^k\|^2_{\Gamma^{-1}-\Lambda_F}
		\)
		holds (surely) for \(k\in\N\) in light of \Cref{thm:geq}.
		The claim then follows by invoking \Cref{thm:decrease} and assert \ref{thm:res}.
	\item\ref{thm:cluster}~
		In the rest of the proof, for conciseness the ``almost sure'' nature of the results will be implied without mention.
		It follows from assert \ref{thm:res} that a subsequence \(\seq{\bm x^k}[k\in K]\) converges to some point \(\bm x^\star\) iff so does the subsequence \(\seq{\bm z^k}[k\in K]\).
		Since \(\T(\bm x^k)\ni\bm z^k\) and both \(\bm x^k\) and \(\bm z^k\) converge to \(\bm x^\star\) as \(K\ni k\to\infty\), the inclusion \(0\in\hat\partial\Phi(\bm x^\star)\) follows from \Cref{thm:critical}.
		Since the full sequences \(\seq{\FBE(\bm x^k)}\) and \(\seq{\Phi(\bm z^k)}\) converge to the same value \(\Phi_\star\) (cf. \Cref{thm:decrease} and assert \ref{thm:decreasez}), due to continuity of \(\FBE\) (\Cref{thm:FBEineq}) it holds that \(\FBE(\bm x^\star)=\Phi_\star\), and in turn the bounds in \Cref{thm:geq} together with assert \ref{thm:res} ensure that \(\Phi(\bm x^\star)=\Phi_\star\) too.
	\qedhere
	\end{proofitemize}
		\end{proof}
	\end{thm}
	
	When \(G\) is convex and \(F\) is strongly convex (that is, each of the functions \(f_i\) is strongly convex), the FBE decreases \(Q\)-linearly in expectation along the iterates generated by the randomized BC-\Cref{alg:BC}.
	
	\begin{thm}[randomized sampling: linear convergence under strong convexity]\label{thm:random:linear}%
		Additionally to \Cref{ass:basic,ass:random}, suppose that \(G\) is convex and that each \(f_i\) is \(\mu_{f_i}\)-strongly convex.
		Then, for all \(k\) the following hold for the iterates generated by \Cref{alg:BC}:
		\begin{subequations}\label{subeq:random:linear}
			\begin{align}
				\label{eq:random:Qlinear}
				\E{\FBE(\bm x^{k+1})-\min\Phi}
			{}\leq{} &
				(1-c)
				\bigl(\FBE(\bm x^k)-\min\Phi\bigr)
			\\
				\E[]{\Phi(\bm z^k)-\min\Phi}
			{}\leq{} &
				\bigl(\Phi(\bm x^0)-\min\Phi\bigr)(1-c)^k
			\\
				\tfrac12\E[]{\|\bm z^k-\bm x^\star\|^2_{\mu_F}}
			{}\leq{} &
				\bigl(\Phi(\bm x^0)-\min\Phi\bigr)(1-c)^k
			\end{align}
		\end{subequations}
		where \(\bm x^\star\coloneqq\argmin\Phi\),
		\(
			\mu_F
		{}\coloneqq{} \tfrac1N
			\blockdiag\bigl(\mu_{f_1}\I_{n_1},\dots\mu_{f_n}\I_{n_N}\bigr)
		\),
		and denoting \(\xi_i=\frac{N-\gamma_iL_{f_i}}{N}\), \(i\in[N]\),
		\begin{equation}\label{eq:cwc}
			c
		{}={}
			\min_{i\in[N]}{
				\set{\tfrac{\xi_ip_i}{\gamma_i}}
				{}\bigg/{}
				\max_{i\in[N]}\set{\tfrac{N-\gamma_i\mu_{f_i}}{\gamma_i^2\mu_{f_i}}}
			}.
		\end{equation}
		Moreover, by setting the stepsizes \(\gamma_i\) and minimum sampling probabilities \(p_i\) as
		\begin{equation}\label{eq:gammaLinear}
			\gamma_i
		{}={}
			\tfrac{N}{\mu_{f_i}}
			\left(1-\sqrt{1-1/\kappa_{i}}\right)
		\quad\text{and}\quad
			p_i
		{}={}
			\frac{
				\left(\sqrt{\kappa_i}+\sqrt{\kappa_i-1}\right)^2
			}{
				\sum_{j=1}^N\left(\sqrt{\kappa_j}+\sqrt{\kappa_j-1}\right)^2
			}
		\end{equation}
		with
		\(
			\kappa_i\coloneqq\frac{L_{f_i}}{\mu_{f_i}}
		\),
		\(i\in[N]\), then the constant \(c\) in \eqref{subeq:random:linear} can be tightened to
		\begin{equation}\label{eq:cbc}
			c
		{}={}
			\tfrac{1}{
				\sum_{i=1}^N\left(
					\sqrt{\kappa_i}+\sqrt{\kappa_i-1}
				\right)^2
			}.
		\end{equation}
		\begin{proof}
	Since \(\bm z^k\) is a minimizer in \eqref{eq:FBE}, the necessary stationarity condition reads
	\(
		\Gamma^{-1}(\bm x^k-\bm z^k)-\nabla F(\bm x^k)
	{}\in{}
		\partial G(\bm z^k)
	\).
	Convexity of \(G\) then implies
	\[
		G(\bm x^\star)
	{}\geq{}
		G(\bm z^k)
		{}+{}
		\innprod{\Gamma^{-1}(\bm x^k-\bm z^k)-\nabla F(\bm x^k)}{\bm x^\star-\bm z^k},
	\]
	whereas from strong convexity of \(F\) we have
	\[
		F(\bm x^\star)
	{}\geq{}
		F(\bm x^k)
		{}+{}
		\innprod{\nabla F(\bm x^k)}{\bm x^\star-\bm x^k}
		{}+{}
		\tfrac12\|\bm x^k-\bm x^\star\|^2_{\mu_F}.
	\]
	By combining these inequalities into \eqref{eq:FBEz}, and denoting \(\Phi_\star\coloneqq\min\Phi=\min\FBE\) (cf. \Cref{thm:min}), we have
	\begin{align*}
		\FBE(\bm x^k)-\Phi_\star
	{}\leq{} &
		\tfrac12\|\bm z^k-\bm x^k\|^2_{\Gamma^{-1}}
		{}-{}
		\tfrac12\|\bm x^\star-\bm x^k\|^2_{\mu_F}
		{}+{}
		\innprod{\Gamma^{-1}(\bm z^k-\bm x^k)}{\bm x^\star-\bm z^k}
	\\
	{}={} &
		\tfrac12\|\bm z^k-\bm x^k\|_{\Gamma^{-1}-\mu_F}^2
		{}+{}
		\innprod{(\Gamma^{-1}-\mu_F)(\bm z^k-\bm x^k)}{\bm x^\star-\bm z^k}
		{}-{}
		\tfrac12\|\bm x^\star-\bm z^k\|_{\mu_F}^2.
	\end{align*}
	Next, by using the inequality
	\(
		\innprod{\bm a}{\bm b}
	{}\leq{}
		\tfrac12\|\bm a\|_{\mu_F}^2
		{}+{}
		\tfrac12\|\bm b\|^2_{\mu_F^{-1}}
	\)
	to cancel out the last term, we obtain
	\begin{align*}
		\FBE(\bm x^k)-\Phi_\star
	{}\leq{} &
		\tfrac12\|\bm z^k-\bm x^k\|_{\Gamma^{-1}-\mu_F}^2
		{}+{}
		\tfrac12\|(\Gamma^{-1}-\mu_F)(\bm x^k-\bm z^k)\|_{\mu_F^{-1}}^2
	\\
	{}={} &
		\tfrac12\|\bm z^k-\bm x^k\|_{\Gamma^{-2}\mu_F^{-1}(\I-\Gamma\mu_F)}^2,
	\numberthis\label{eq:QUB}
	\end{align*}
	where the last identity uses the fact that the matrices are diagonal.
	Combined with \eqref{eq:EFBE+} the claimed \(Q\)-linear convergence \eqref{eq:random:Qlinear} with factor \(c\) as in \eqref{eq:cwc} is obtained.
	The $R$-linear rates in terms of the cost function and distance from the solution are obtained by repeated application of \eqref{eq:random:Qlinear} after taking (unconditional) expectation from both sides and using \Cref{thm:FBEineq}. 
	
	To obtain the tighter estimate \eqref{eq:cbc}, observe that \eqref{eq:EFBE+} with the choice
	\[
	\textstyle
		p_i
	{}\coloneqq{}
		\tfrac{1}{\gamma_i\mu_{f_i}}
		\tfrac{N-\gamma_i\mu_{f_i}}{N-\gamma_iL_{f_i}}
		\left(
			\sum_j{
				\tfrac{1}{\gamma_j\mu_{f_j}}
				\tfrac{N-\gamma_j\mu_{f_j}}{N-\gamma_jL_{f_j}}
			}
		\right)^{-1},
	\]
	which equals the one in \eqref{eq:gammaLinear} with $\gamma_i$ as prescribed, yields
	\begin{align*}
		\E{\FBE(\bm x^{k+1})-\Phi_\star}
	{}\leq{} &
		\FBE(\bm x^k)-\Phi_\star
		{}-{}
		\left(
			\textstyle
			2N\sum_j{
				\tfrac{1}{\gamma_j\mu_j}
				\tfrac{N-\gamma_j\mu_j}{N-\gamma_jL_j}
			}
		\right)^{-1}
		\sum_{i=1}^N{
			\tfrac{N-\gamma_i\mu_{f_i}}{\gamma_i^2\mu_{f_i}}
			\|z_i^k-x_i^k\|^2
		}
	\\
	{}={} &
		\textstyle
		\FBE(\bm x^k)-\Phi_\star
		{}-{}
		\left(
			2N\sum_j{
				\tfrac{1}{\gamma_j\mu_j}
				\tfrac{N-\gamma_j\mu_j}{N-\gamma_jL_j}
			}
		\right)^{-1}
		\|\bm z^k-\bm x^k\|_{\Gamma^{-1}\mu_F^{-1}(\Gamma^{-1}-\mu_F)}^2.
	\end{align*}
	The assert now follows by combining this with \eqref{eq:QUB} and replacing the values of \(\gamma_i\) as proposed in \eqref{eq:gammaLinear}.
		\end{proof}
	\end{thm}
	Notice that as \(\kappa_i\)'s approach \(1\) the linear rate tends to \(1-\nicefrac1N\).

		\subsection{Cyclic, shuffled and essentially cyclic samplings}\label{sec:cyclic}

	
	In this section we analyze the convergence of the BC-\Cref{alg:BC} when a cyclic, shuffled cyclic or (more generally) an essentially cyclic sampling \cite{tseng1987relaxation,tseng2001convergence,hong2017iteration,chow2017cyclic,xu2017globally} is used.
	As formalized in the following standing assumption, an additional convexity requirement for the nonsmooth term \(G\) is needed.
	
	\begin{ass}[essentially cyclic sampling requirements]\label{ass:cyclic}%
		In problem \eqref{eq:P}, function \(G\) is convex.
		Moreover, there exists $T\geq 1$ such that in \Cref{alg:BC} each index is selected at least once within any interval of $T$ iterations.
	\end{ass}
	
	Note that having \(T<N\) is possible because of our general sampling strategy where sets of indices can be sampled within the same iteration.
	For instance, \(T=1\) corresponds to \(I^{k+1}=[N]\) for all \(k\), in which case \Cref{alg:BC} would reduce to a (full) proximal gradient scheme.
	
	Two notable special cases of single index selection rules are the cyclic and shuffled cyclic sampling strategies. 
	\begin{itemize}[%
		leftmargin=*,
		label={},
		itemindent=0cm,
		labelsep=0pt,
		partopsep=0pt,
		parsep=0pt,
		listparindent=0pt,
		topsep=0pt,
	]
	\item{\sc Shuffled cyclic sampling:} corresponds to setting 
		\begin{equation}\label{eq:ShufCyclicRule}
			I^{k+1}=\set{\pi_{\lfloor\nicefrac kN\rfloor}\bigl(\mod(k,N)+1\bigr)}\quad \text{for all}\quad k\in\N,
		\end{equation}
		where $\pi_0,\pi_1,\dots$ are permutations of the set of indices $[N]$ (chosen randomly or deterministically).
	\item{\sc Cyclic sampling:} corresponds to the case \eqref{eq:ShufCyclicRule} with $\pi_{\lfloor\nicefrac kN\rfloor}=\id$, \ie,   
		\begin{equation}\label{eq:cyclicRule}
			I^{k+1}=\set{\mod(k,N)+1}\quad \text{for all}\quad k\in\N.
		\end{equation}
	\end{itemize}
	We remark that in practice it has been observed that an effective sampling technique is to use random shuffling after each cycle \cite[\S2]{bertsekas2015convex}. 
	Consistently with the deterministic nature of the essentially cyclic sampling, all results of the previous section hold surely, as opposed to almost surely.
	
	\begin{thm}[essentially cyclic sampling: subsequential convergence]\label{thm:cyclic:subseq}%
		Suppose that \Cref{ass:basic,ass:cyclic} are satisfied.
		Then, all the asserts of \Cref{thm:random:subseq} hold surely.
		\begin{proof}
	We first establish an important descent inequality for \(\FBE\) after every $T$ iterations, cf. \eqref{eq:Essential_cyclic_descent}.
	Convexity of \(G\), entailing $\prox_{G}^{\Gamma^{-1}}$ being Lipschitz continuous (cf. \Cref{thm:FNE}), allows the employment of techniques similar to those in \cite[Lemma 3.3]{beck2013convergence}. 
	Since all indices are updated at least once every \(T\) iterations, one has that
	\begin{equation}\label{eq:ki}
		\ki
	{}\coloneqq{}
		\min\set{t\in[T]}[
			\text{\(i\) is sampled at iteration \(T\nu+t-1\)}
		]
	\end{equation}
	is well defined for each index \(i\in[N]\) and \(\nu\in\N\).
	Since \(i\) is sampled at iteration \(T\nu+\ki-1\) and \(x_i^{T\nu}=x_i^{T\nu+1}=\dots=x_i^{T\nu+\ki-1}\) by definition of \(\ki\), it holds that
	\begin{align*}
		x_i^{T\nu+\ki}
	{}={} &
		x_i^{T\nu+\ki-1}
		{}+{}
		\trans{U_i}\,\left(
			\T(\bm x^{T\nu+\ki-1})
			{}-{}
			\bm x^{T\nu+\ki-1}
		\right)
	\\
	\numberthis\label{eq:equiiter_ki}
	{}={} &
		x_i^{T\nu}
		{}+{}
		\trans{U_i}\,\left(
			\T(\bm x^{T\nu+\ki-1})
			{}-{}
			\bm x^{T\nu+\ki-1}
		\right),
	\end{align*}
	\ifaccel\else
	where $U_i\in \R^{(\sum_jn_j)\times n_i}$ denotes the $i$-th block column of the identity matrix so that for a vector $v\in \R^{n_i}$
	\begin{equation}\label{eq:U}
		U_iv
	{}={}
		\trans{(0,\dots, 0, \!\overbracket{\,v\,}^{\mathclap{i\text{-th}}}\!, 0, \dots, 0)}.
	\end{equation}
	\fi
	For all $t\in[T]$ the following holds
	\begin{align*}
		\FBE(\bm x^{T(\nu+1)})
		{}-{}
		\FBE(\bm x^{T\nu})
	{}={} &
		\sum_{\tau=1}^T\left(
			\FBE(\bm x^{T\nu+\tau})
			{}-{}
			\FBE(\bm x^{T\nu+\tau-1})
		\right)
	\\
	{}\leq{} &
		\FBE(\bm x^{T\nu+t})
		{}-{}
		\FBE(\bm x^{T\nu+t-1})
	\\
	{}\leq{} &
		-\tfrac{\xi_{\rm min}}{2}
		\|\bm x^{T\nu+t}-\bm x^{T\nu+t-1}\|_{\Gamma^{-1}}^2,
	\numberthis\label{eq:descent_esscyc}
	\end{align*}
	where $\xi_i\coloneqq \tfrac{N-\gamma_{i}L_{f_{i}}}{N}$ as in \Cref{thm:Igeq}, $\xi_{\rm min}\coloneqq \min_{i\in[N]}\set{\xi_i}$, and the two inequalities follow from \Cref{thm:Igeq}.
	Moreover, using triangular inequality for $i\in[N]$ yields
	\begin{align*}
		\|\bm x^{T\nu+\ki-1}-\bm x^{T\nu}\|_{\Gamma^{-1}}
	{}\leq{} &
		\sum_{\tau=1}^{\ki-1}\|\bm x^{T\nu+\tau}-\bm x^{T\nu+\tau-1}\|_{\Gamma^{-1}}
	\\
	\numberthis\label{eq:new25}
	{}\leq{} &
		\tfrac{T}{\sqrt{\xi_{\rm min}\nicefrac{}{2}}}
		\left(
			\FBE(\bm x^{T\nu})
			{}-{}
			\FBE(\bm x^{T(\nu+1)})
		\right)^{\nicefrac12},
	\end{align*}
	where the second inequality follows from \eqref{eq:descent_esscyc} together with the fact that \(\ki\leq T\).
	For all \(i\in[N]\), from the triangular inequality and the \(L_{\bf T}\)-Lipschitz continuity of \(\T\) (\Cref{thm:TLip}) we have
	\begin{align*}
		\gamma_i^{-\nicefrac12}
		\|\trans{U_i}\,(\bm x^{T\nu}-\T(\bm x^{T\nu}))\|
	{}\leq{} &
		\gamma_i^{-\nicefrac12}
		\|\trans{U_i}\,\bigl(\bm x^{T\nu}-\T(\bm x^{T\nu+\ki-1})\bigr)\|
	\\
	&
		{}+{}
		\gamma_i^{-\nicefrac12}
		\|\trans{U_i}\,\bigl(\T(\bm x^{T\nu+\ki-1})-\T(\bm x^{T\nu})\bigr)\|
	\\
	{}\leq{} &
		\gamma_i^{-\nicefrac12}
		\|x_i^{T\nu+\ki-1}-x_i^{T\nu+\ki}\|
	\\
	&
		{}+{}
		\|\T(\bm x^{T\nu+\ki-1})-\T(\bm x^{T\nu})\|_{\Gamma^{-1}}
	\\
	{}\leq{} &
		\|\bm x^{T\nu+\ki-1}-\bm x^{T\nu+\ki}\|_{\Gamma^{-1}}
		{}+{}
		L_{\bf T}\|\bm x^{T\nu+\ki-1}-\bm x^{T\nu}\|_{\Gamma^{-1}}
	\\
	\numberthis\label{eq:sqrtbound}
	{}\overrel[\leq]{\eqref{eq:descent_esscyc},~\eqref{eq:new25}}{} &
		\tfrac{1+TL_{\bf T}}{\sqrt{\xi_{\rm min}/2}}
		\left(
			\FBE(\bm x^{T\nu})
			{}-{}
			\FBE(\bm x^{T(\nu+1)})
		\right)^{\nicefrac12}.
	\end{align*}
	By squaring and summing over \(i\in[N]\), we obtain
	\begin{equation}\label{eq:Essential_cyclic_descent}
		\FBE(\bm x^{T(\nu+1)})-\FBE(\bm x^{T\nu})
	{}\leq{}
		-\tfrac{\xi_{\rm min}}{2N(1+TL_{\bf T})^2}
		\|\bm z^{T\nu}-\bm x^{T\nu}\|^{2}_{\Gamma^{-1}}.
	\end{equation}
	By telescoping the inequality and using the fact that \(\min\FBE=\min\Phi\)
	\ifarxiv
		by
	\else
		shown in
	\fi
	\Cref{thm:min}, we obtain that
	\(
		\seq{\|\bm z^{T\nu}-\bm x^{T\nu}\|^2_{\Gamma^{-1}}}[\nu\in\N]
	\)
	has finite sum, and in particular vanishes.
	Clearly, by suitably shifting, for every \(t\in[T]\) the same can be said for the sequence
	\(
		\seq{\|\bm z^{T\nu+t}-\bm x^{T\nu+t}\|^2_{\Gamma^{-1}}}[\nu\in\N]
	\).
	The whole sequence
	\(
		\seq{\|\bm z^k-\bm x^k\|^2}
	\)
	is thus summable, and we may now infer the claim as done in the proof of \Cref{thm:random:subseq}.
		\end{proof}
	\end{thm}
	
	In the next theorem explicit linear convergence rates are derived under the additional strong convexity assumption for the smooth functions.
	The cyclic and shuffled cyclic cases are treated separately, as tighter bounds can be obtained 	by leveraging the fact that within cycles of \(N\) iterations every index is updated exactly once.
	
	\begin{thm}[essentially cyclic sampling: linear convergence under strong convexity]\label{thm:cyclic:linear}%
		Additionally to \Cref{ass:basic,ass:cyclic}, suppose that each function \(f_i\) is \(\mu_{f_i}\)-strongly convex.
		Then, denoting
		\(
			\delta
		{}\coloneqq{}
			\min_{i\in[N]}\set{
				\tfrac{\gamma_i\mu_{f_i}}{N}
			}
		\)
		and
		\(
			\Delta
		{}\coloneqq{}
			\max_{i\in[N]}\set{
				\tfrac{\gamma_iL_{f_i}}{N}
			}
		\),
		for all \(\nu\in\N\) the following hold for the iterates generated by \Cref{alg:BC}:
		\begin{subequations}\label{subeq:cyclic:linear}
			\begin{align}
				\label{eq:cyclic:Qlinear}
				\FBE(\bm x^{T(\nu+1)})-\min\Phi
			{}\leq{} &
				(1-c)
				\bigl(\FBE(\bm x^{T\nu})-\min\Phi\bigr)
			\\ \label{eq:cyclic:Rlinear1}
				\Phi(\bm z^{T\nu})-\min\Phi
			{}\leq{} &
				\bigl(\Phi(\bm x^0)-\min\Phi\bigr)(1-c)^\nu
			\\ \label{eq:cyclic:Rlinear2}
				\tfrac12\|\bm z^{T\nu}-\bm x^\star\|^2_{\mu_F}
			{}\leq{} &
				\bigl(\Phi(\bm x^0)-\min\Phi\bigr)(1-c)^\nu
			\end{align}
		\end{subequations}
		where \(\bm x^\star\coloneqq\argmin\Phi\),
		\(
			\mu_F
		{}\coloneqq{} \tfrac1N
			\blockdiag\bigl(\mu_{f_1}\I_{n_1},\dots\mu_{f_n}\I_{n_N}\bigr)
		\),
		and
		\begin{equation}\label{eq:cyclic:cwc}
				c
			{}={}
				\frac{
					\delta(1-\Delta)
				}{
					N\bigl(1+T(1-\delta)\bigr)^2
					(1-\delta)
				}.
		\end{equation}
		In the case of shuffled cyclic \eqref{eq:ShufCyclicRule} or cyclic \eqref{eq:cyclicRule} sampling, the inequalities can be tightened by
		replacing \(T\) with \(N\) and with
		\begin{equation}\label{eq:linearShuffledCyclic}
			c
		{}={}
			\frac{\delta(1-\Delta)}{N\left(2-\delta\right)^{2}\left(1-\delta\right)}.
		\end{equation}
		\begin{proof}
	\begin{proofitemize}
	\item\emph{The general essentially cyclic case.}~
		Since \(\T\) is \(L_{\bf T}\)-Lipschitz continuous with \(L_{\bf T}=1-\delta\) as shown in \Cref{thm:contractive}, inequality \eqref{eq:Essential_cyclic_descent} becomes 
		\[
			\FBE(\bm x^{T(\nu+1)})
			{}-{}
			\FBE(\bm x^{T\nu})
		{}\leq{}
			-\tfrac{1-\Delta}{2N(1+T(1-\delta))^2}
			\|\bm z^{T\nu}-\bm x^{T\nu}\|^2_{\Gamma^{-1}}.
		\]
		Moreover, it follows from \eqref{eq:QUB} that
		\begin{equation}\label{eq:strongLB}
			\FBE(\bm x^{T\nu})-\Phi_\star
		{}\leq{}
			\tfrac12
			(\delta^{-1}-1)
			\|\bm z^{T\nu}-\bm x^{T\nu}\|_{\Gamma^{-1}}^2.
		\end{equation}
		By combining the two inequalities the claimed \(Q\)-linear convergence \eqref{eq:cyclic:Qlinear} with factor \(c\) as in \eqref{eq:cyclic:cwc} is obtained.
		In turn, the $R$-linear rates \eqref{eq:cyclic:Rlinear1} and \eqref{eq:cyclic:Rlinear2} follow from \Cref{thm:FBEineq}. 
	\item\emph{The shuffled cyclic case.}~
		Let us now suppose that the sampling strategy follows a shuffled rule as in \eqref{eq:ShufCyclicRule} with permutations \(\pi_0,\pi_1,\dots\) (hence in the cyclic case $\pi_\nu=\id$ for all $\nu\in\N$).
		Let $U_i$ be as in \eqref{eq:U} and $\xi_{\rm min}$ as in the proof of \Cref{thm:cyclic:subseq}.  Observe that \(\ki=\pi_\nu^{-1}(i)\leq N\) for \(\ki\) as defined in \eqref{eq:ki}.
		For all $t\in[N]$
		\begin{align*}
			\FBE(\bm x^{N(\nu+1)}) -\FBE(\bm x^{N\nu})
		{}\leq{} &
			\FBE(\bm x^{N\nu+t-1}) -\FBE(\bm x^{N\nu})
		\\
		{}\leq{} &
			-\tfrac{\xi_{\rm min}}{2}\sum_{\tau=1}^{t-1} \|\bm x^{N\nu+\tau}-\bm x^{N\nu+\tau-1}\|^2_{\Gamma^{-1}} 
		\\
		\numberthis\label{eq:tighterNoTri}
		{}={} &
			 -\tfrac{\xi_{\rm min}}{2}\|\bm x^{N\nu+t-1}-\bm x^{N\nu}\|^2_{\Gamma^{-1}},
		\end{align*}
		where the equality follows from the fact that at every iteration a different coordinate is updated (and that $\Gamma$ is diagonal), and the inequalities from \Cref{thm:Igeq}. Similarly, \eqref{eq:descent_esscyc} holds with $T$ replaced by \(N\) (despite the fact that \(T\) is not necessarily \(N\), but is rather bounded as \(T\leq 2N-1\)).
		By using \eqref{eq:tighterNoTri} in place of \eqref{eq:new25}, inequality \eqref{eq:sqrtbound} is tightened as follows 
		\[
			\gamma_i^{-\nicefrac12}
			\|\trans{U_i}(\bm x^{N\nu}-\T(\bm x^{N\nu}))\|
		{}\leq{}
			\tfrac{1+L_{\bf T}}{\sqrt{\xi_{\rm min}/2}}
			\left(
				\FBE(\bm x^{N\nu})
				{}-{}
				\FBE(\bm x^{N(\nu+1)})
			\right)^{\nicefrac12}.
		\]
		By squaring and summing for \(i\in[N]\) we obtain
		\begin{equation}\label{eq:cyclic_descent}
			\FBE(\bm x^{N(\nu+1)})-\FBE(\bm x^{N\nu})
		{}\leq{}
			-\tfrac{\xi_{\rm min}}{2N(1+L_{\bf T})^2}
			\|\bm z^{N\nu}-\bm x^{N\nu}\|^2_{\Gamma^{-1}}
		{}={}
			-\tfrac{1-\Delta}{2N(1+L_{\bf T})^2}
			\|\bm z^{N\nu}-\bm x^{N\nu}\|^2_{\Gamma^{-1}},
		\end{equation}
		where \(L_{\bf T}=1-\delta\) as discussed above.
		By combining this and \eqref{eq:strongLB} (with \(T\) replaced by \(N\)) the improved coefficient \eqref{eq:linearShuffledCyclic} is obtained.
	\qedhere
	\end{proofitemize}
		\end{proof}
	\end{thm}
	
	Note that if one sets $\gamma_i = \alpha N/L_{f_i}$ for some $\alpha\in(0,1)$, then  $\delta = \alpha\min_{i\in[N]} \set{\nicefrac{\mu_{f_i}}{L_{f_i}}}$ and $\Delta=\alpha$.
	With this selection, as the condition number approaches $1$ the rate in \eqref{eq:linearShuffledCyclic} tends to $1-\frac{\alpha}{N\left(2-\alpha\right)^{2}}$.   

		\subsection{Global and linear convergence with KL inequality}

	The convergence analyses of the randomized and essentially cyclic cases both rely on a descent property on the FBE that quantifies the progress in the minization of \(\FBE\) in terms of the squared forward-backward residual \(\|\bm x-\bm z\|^2\).
	A subtle but important difference, however, is that the inequality \eqref{eq:ExSD} in the former case involves a conditional expectation, whereas \eqref{eq:Essential_cyclic_descent} in the latter does not.
	The \emph{sure} descent property occurring for essentially cyclic sampling strategies is the key for establishing global (as opposed to subsequential) convergence based on the Kurdyka-\L ojasiewicz (KL) property \cite{lojasiewicz1963propriete,lojasiewicz1993geometrie,kurdyka1998gradients}.
	A similar result is achieved in \cite{xu2017globally}, which however considers the complementary case to problem \eqref{eq:P} where the nonsmooth function \(G\) is assumed to be separable, and thus the cost function itself can serve as Lyapunov function.
	
	\begin{defin}[KL property with exponent \(\theta\)]\label{def:KL}%
		A proper lsc function \(\func{h}{\R^n}{\Rinf}\) is said to have the \DEF{Kurdyka-{\L}ojasiewicz} (KL) property with exponent \(\theta\in(0,1)\) at \(\bar w\in\dom h\) if there exist \(\varepsilon,\eta,\varrho>0\) such that
		\[
			\psi'(h(w)-h(\bar w))\dist(0,\partial h(w))\geq 1
		\]
		holds for all \(w\) such that \(\|w-\bar w\|<\varepsilon\) and \(h(\bar w)<h(w)<h(\bar w)+\eta\), where \(\psi(s)\coloneqq\varrho s^{1-\theta}\).
		We say that \(h\) satisfies the KL property with exponent \(\theta\) (without mention of \(\bar w\)) if it satisfies the KL property with exponent \(\theta\) at any \(\bar w\in\dom\partial h\).
	\end{defin}
	
	Semialgebraic functions comprise a wide class of functions that enjoy this property \cite{bolte2007clarke,bolte2007lojasiewicz}, which has been extensively exploited to provide convergence rates of optimization algorithms \cite{attouch2009convergence,attouch2010proximal,attouch2013convergence,bolte2014proximal,frankel2015splitting,ochs2014ipiano,li2016douglas,xu2013block}.
	Based on this, in the next result we provide sufficient conditions ensuring global and \(R\)-linear convergence of \Cref{alg:BC} with essentially cyclic sampling.
	
	\begin{thm}[essentially cyclic sampling: global and linear convergence]\label{thm:cyclic:global}%
		Additionally to \Cref{ass:basic,ass:cyclic}, suppose that \(\Phi\) has the KL property with exponent \(\theta\in(0,1)\) (as is the case when \(f_i\) and \(G\) are semialgebraic), and is coercive.
		Then, any sequences \(\seq{\bm x^k}\) and \(\seq{\bm z^k}\) generated by \Cref{alg:BC} converge to (the same) stationary point \(\bm x^\star\).
		Moreover, if \(\theta\leq\nicefrac12\) then \(\seq{\|\bm z^k-\bm x^k\|}\), \(\seq{\bm x^k}\) and \(\seq{\bm z^k}\) converge at $R$-linear rate.
		\begin{proof}
	Let \(\seq{\bm x^k}\) and \(\seq{\bm z^k}\) be sequences generated by \Cref{alg:BC} with essentially cyclic sampling, and let \(\Phi_\star\) be the limit of the sequence \(\seq{\FBE(\bm x^k)}\) as in \Cref{thm:decrease}.
	To avoid trivialities, we may assume that \(\FBE(\bm x^k)\gneqq\Phi_\star\) for all \(k\), for otherwise the sequence \(\seq{\bm x^k}\) is asymptotically constant, and thus so is $\seq{\bm z^k}$.
	Let \(\Omega\) be the set of accumulation points of \(\seq{\bm x^k}\), which is compact and such that \(\FBE\equiv\Phi_\star\) on \(\Omega\), as ensured by \Cref{thm:cyclic:subseq}.
	It follows from \Cref{thm:loja} and \cite[Lem. 1(ii)]{attouch2009convergence} that \(\FBE\) enjoys a \emph{uniform} KL property on \(\Omega\); in particular,
	\(
		\psi'(\FBE(\bm x^k)-\Phi_\star)\dist(0,\partial\FBE(\bm x^k))
	{}\geq{}
		1
	\)
	holds for all \(k\) large enough such that \(\bm x^k\) is sufficiently close to \(\Omega\) and \(\FBE(\bm x^k)\) is sufficiently close to \(\Phi_\star\), where \(\psi(s)=\varrho s^{1-\theta'}\) for some \(\varrho>0\) and \(\theta'=\max\set{\theta,\nicefrac12}\).
	Combined with \Cref{thm:subdiffdist}, for all \(k\) large enough we thus have
	\begin{equation}\label{eq:KL}
		\psi'(\FBE(\bm x^k)-\Phi_\star)
	{}\geq{}
		\frac{c}{\|\bm x^k-\bm z^k\|_{\Gamma^{-1}}},
	\end{equation}
	where
	\(
		c
	{}\coloneqq{}
		\frac{
			N\min_i\set{\sqrt{\gamma_i}}
		}{
			N+\max_i\set{\gamma_iL_{f_i}}
		}
	{}>{}
		0
	\).
	Let
	\(
		\Delta_k\coloneqq\psi(\FBE(\bm x^k)-\Phi_\star)
	\).
	By combining \eqref{eq:KL} and \eqref{eq:Essential_cyclic_descent} we have that there exists a constant \(c'>0\) such that
	\begin{equation}\label{eq:KLinequality}
		\Delta_{(\nu+1)T}
		{}-{}
		\Delta_{\nu T}
	{}\leq{}
		\psi'(\FBE(\bm x^{\nu T})-\Phi_\star)
		\left(\FBE(\bm x^{(\nu+1)T})-\FBE(\bm x^{\nu T})\right)
	{}\leq{}
		-c'
		\|\bm x^{\nu T}-\bm z^{\nu T}\|_{\Gamma^{-1}}
	\end{equation}
	holds for all $\nu\in\N$ large enough (the first inequality uses concavity of \(\psi\)).
	By summing over \(\nu\) (sure) summability of the sequence \(\seq{\|\bm x^{\nu T}-\bm z^{\nu T}\|}[\nu\in\N]\) is obtained.
	By suitably shifting, for every \(t\in[T]\) the same can be said for the sequence
	\(
		\seq{\|\bm z^{T\nu+t}-\bm x^{T\nu+t}\|}[\nu\in\N]
	\),
	and since \(T\) is finite we conclude that the whole sequence
	\(
		\seq{\|\bm z^k-\bm x^k\|}
	\)
	is summable.
	Since \(\|\bm x^{k+1}-\bm x^k\|\leq\|\bm z^k-\bm x^k\|\) we conclude that \(\seq{\bm x^k}\) has finite length and is thus convergent (to a single point), and consequently so is \(\seq{\bm z^k}\).
	Suppose now that \(\theta\leq\nicefrac12\), so that \(\psi(s)=\varrho\sqrt s\).
	Then,
	\[
		\|\bm x^{\nu T}-\bm z^{\nu T}\|_{\Gamma^{-1}}
	{}\overrel[\geq]{\eqref{eq:KL}}{}
		\tfrac{2c}{\varrho}
		\sqrt{\FBE(\bm x^{\nu T})-\Phi_\star}
	{}={}
		\tfrac{2c}{\varrho^2}
		\psi(\FBE(\bm x^{\nu T})-\Phi_\star)
	{}={}
		\tfrac{2c}{\varrho^2}
		\Delta_{\nu T}.
	\]
	Combined with \eqref{eq:KLinequality} it follows that $\seq{\Delta_{\nu T}}[\nu\in\N]$ conveges $Q$-linearly.
	By rearranging \eqref{eq:KLinequality} as
	\[
		c'\|\bm x^{\nu T}-\bm z^{\nu T}\|_{\Gamma^{-1}}
	{}\leq{}
		\Delta_{\nu T}
		{}-{}
		\Delta_{(\nu+1)T}
	{}\leq{}
		\Delta_{\nu T},
	\]
	\(R\)-linear convergence of \(\seq{\|\bm x^{\nu T}-\bm z^{\nu T}\|}[\nu\in\N]\) follows.
	By suitably shifting, for every \(t\in[T]\) the same can be said for the sequence
	\(
		\seq{\|\bm z^{T\nu+t}-\bm x^{T\nu+t}\|}[\nu\in\N]
	\),
	and since \(T\) is finite we conclude that the whole sequence
	\(
		\seq{\|\bm z^k-\bm x^k\|}
	\)
	converges $R$-linearly. On the other hand, since \(\|\bm x^{k+1}-\bm x^k\|\leq\|\bm z^k-\bm x^k\|\), also \(\seq{\|\bm x^{k+1}-\bm x^k\|}\) converges \(R\)-linearly, hence so does \(\seq{\bm x^k}\). By combining the two, we conclude that also $\seq{\bm z^k}$ converges $R$-linearly. 
		\end{proof}
	\end{thm}

	\section{Nonconvex finite sum problems: the Finito/MISO algorithm}\label{sec:Finito}

	As mentioned in \Cref{sec:Introduction}, if \(G\) is of the form \eqref{eq:FINITOG} then problem \eqref{eq:P} reduces to the finite sum minimization presented in \eqref{eq:FSP}.
	Most importantly, the proximal mapping of the original nonsmooth function \(G\) can be easily expressed in terms of that of the small function \(g\) in the reduced finite sum reformulation, as shown in the next lemma.
	
	\begin{lem}
		Given \(\gamma_i>0\), \(i\in[N]\), let
		\(
			\Gamma
		{}\coloneqq{}
			\blockdiag(\gamma_1I_n,\dots,\gamma_NI_n)
		\)
		and
		\(
			\hat\gamma
		{}\coloneqq{}
			\bigl(\sum_{i=1}^N\gamma_i^{-1}\bigr)^{-1}
		\).
		Then, for \(G\) as in \eqref{eq:FINITOG} and any \(\bm u\in\R^{Nn}\)
		\[
		\textstyle
			\prox_G^{\Gamma^{-1}}(\bm u)
		{}={}
			\set{(\hat v,\dots,\hat v)}[
				\hat v
			{}\in{}
				\prox_{\hat\gamma g}(\hat u)
			]
		\quad\text{where}\quad
			\hat u
		{}\coloneqq{}
			\hat\gamma
			\sum_{i=1}^N\gamma_i^{-1}u_i.
		\]
		\begin{proof}
			Observe first that for every \(w\in\R^n\) one has
			\begin{align*}
			\textstyle
				\sum_i\gamma_i^{-1}\|w-u_i\|^2
			{}={} &
			\textstyle
				\sum_i\gamma_i^{-1}\|\hat u-u_i\|^2
				{}+{}
				\sum_i\gamma_i^{-1}\|w-\hat u\|^2
				{}+{}
				\smashoverbrace{
					\textstyle
					2\sum_i\gamma_i^{-1}\innprod{\hat u-u_i}{w-\hat u}
				}{
					=0
				}
			\\
			\numberthis\label{eq:mean}
			{}={} &
			\textstyle
				\sum_i\gamma_i^{-1}\|\hat u-u_i\|^2
				{}+{}
				\hat\gamma^{-1}\|w-\hat u\|^2.
			\end{align*}
			Next, observe that since \(\dom G\subseteq C\) (the consensus set),
			\begin{align*}
				\prox_G^{\Gamma^{-1}}(\bm u)
			{}={} &
				\argmin_{\bm w\in\R^{Nn}}\set{\textstyle
					G(\bm w)+\sum_{i=1}^N\tfrac{1}{2\gamma_i}\|w_i-u_i\|^2
				}
			\\
			{}={} &
				\argmin_{\bm w\in\R^{Nn}}\set{\textstyle
					G(\bm w)+\sum_{i=1}^N\tfrac{1}{2\gamma_i}\|w_i-u_i\|^2
				}[
					w_1=\dots=w_N
				]
			\\
			{}={} &
				\argmin_{(w,\dots,w)}\set{\textstyle
					g(w)+\sum_{i=1}^N\tfrac{1}{2\gamma_i}\|w-u_i\|^2
				}
			\\
			{}\overrel*{\eqref{eq:mean}}{} &
				\argmin_{(w,\dots,w)}\set{\textstyle
					g(w)
					{}+{}
					\tfrac{1}{2\hat\gamma}\|w-\hat u\|^2
				}
			{}={}
				\set{(\hat v,\dots,\hat v)}[
					\hat v\in\prox_{\hat\gamma g}(\hat u)
				]
			\end{align*}
			as claimed.
		\end{proof}
	\end{lem}
	
	If all stepsizes are set to a same value \(\gamma\), so that \(\Gamma=\gamma\I_{Nn}\), then the forward-backward step reduces to
	\begin{align*}
		\bm z
		{}\in{}
		\prox_G^{\Gamma^{-1}}(\bm x-\Gamma\nabla F(\bm x))
	\quad\Leftrightarrow\quad &
		\bm z=(\bar z,\dots,\bar z),
	\\
	\numberthis\label{eq:FBFinito}
	&
		\bar z
		{}\in{}
		\prox_{\gamma g\nicefrac{}N}\left(
			\textstyle
			\tfrac1N\sum_{j=1}^N\bigl(
				x_j-\tfrac\gamma N\nabla f_j(x_j)
			\bigr)
		\right).
	\end{align*}
	The argument of \(\prox_{\gamma g\nicefrac{}{N}}\) is the (unweighted) average of the forward operator.
	By applying \Cref{alg:BC} with \eqref{eq:FBFinito}, Finito/MISO \cite{defazio2014finito,mairal2015incremental} is recovered.
	Differently from the existing convergence analyses, ours covers fully nonconvex and nonsmooth problems, more general sampling strategies and the possibility to select different stepsizes \(\gamma_i\) for each block, which can have a significant impact on the performance compared to the case where all stepsizes are equal.
	Moreover, to the best of our knowledge this is the first work that shows global convergence and linear rates even when the smooth functions are nonconvex.
	The resulting scheme is presented in \Cref{alg:Finito}.
	We remark that the consensus formulation to recover Finito/MISO (although from a different umbrella algorithm) was also observed in \cite{davis2016smart} in the convex case.
	Moreover, the Finito/MISO algorithm with cyclic sampling is also studied in \cite{mokhtari2018surpassing} when \(g\equiv0\) and $f_i$ are strongly convex functions; consistently with \Cref{ass:cyclic}, our analysis covers the more general essentially cyclic sampling even in the presence of a nonsmooth convex term $g$ and allowing the smooth functions $f_i$ to be nonconvex.
	Randomized Finito/MISO with $g\equiv 0$ is also studied in the recent work \cite{qian2019miso}; although their analysis is limited to a single stepsize, in the convex case it is allowed to be larger than our worst-case stepsize \(\min_i\gamma_i\).
	
	\begin{algorithm}
		\caption{Nonconvex proximal Finito/MISO for problem \eqref{eq:FSP}
		}
		\label{alg:Finito}
	\begin{algorithmic}[1]
	\item[{\sc Require}]
		\(
			x^{\rm init}\in\R^n
		\),~
		\(\gamma_i\in(0,\nicefrac{N}{L_{f_i}})\),
		{\small \(i\in[N]\)}
	\Statex
		\(
			\hat\gamma
		{}\coloneqq{}
			\bigl(\sum_{i=1}^N\gamma_i^{-1}\bigr)^{-1}
		\),~~
		\(
			s_i
		{}={}
			x^{\rm init}-\frac{\gamma_i}{N}\nabla f_i(x^{\rm init})
		\)~
		\(i\in[N]\),~~
		\(
			\hat s
		{}={}
			{\hat\gamma}\sum_{i=1}^N\gamma_i^{-1}s_i
		\)
	\item[{\sc Repeat} until convergence]
		\State
			select a set of indices \(I\subseteq[N]\)
		\State
			\(
				z
			{}\in{}
				\prox_{\hat\gamma g}(\hat s)
			\)
		\For{ \(i\in I\) }
		\State
			\(
				v
			{}\gets{}
				z-\frac{\gamma_i}{N}\nabla f_i(z)
			\)
		\State
			update~~
			\(
				\hat s
			{}\gets{}
				\hat s+\tfrac{\hat\gamma}{\gamma_i}(v-s_i)
			\)
			~~and~~
			\(
				s_i\gets v
			\)
		\EndFor
	\item[{\sc Return} $z$ ]
	\end{algorithmic}
	\end{algorithm}
	
	The convergence results from \Cref{sec:convergence} are immediately translated to this setting by noting that the bold variable ${\bm z}^k$ corresponds to $(z^k,\dots,z^k)$.
	Therefore, $\Phi({\bm z^k})= \varphi(z^k)$ where $\varphi$ is the cost function for the finite sum problem.
	
	\begin{cor}[subsequential convergence of \Cref{alg:Finito}]\label{thm:Finito:convergence}%
		In the finite sum problem \eqref{eq:FSP} suppose that \(\argmin\varphi\) is nonempty, \(g\) is proper and lsc, and each \(f_i\) is \(L_{f_i}\)-Lipschitz differentiable, \(i\in[N]\).
		Then, the following hold almost surely (resp. surely) for the sequence $\seq{z^k}$ generated by \Cref{alg:Finito} with randomized sampling strategy as in \Cref{ass:random} (resp. with any essentially cyclic sampling strategy and $g$ convex as required in \Cref{ass:cyclic}):
		\begin{enumerate}
		\item
			the sequence \(\seq{\varphi(z^k)}\) converges to a finite value \(\varphi_\star\leq\varphi(x^{\rm init})\);
		\item
			all cluster points of the sequence \(\seq{z^k}\) are stationary and on which \(\varphi\) equals \(\varphi_\star\).
		\end{enumerate}
		If, additionally, \(\varphi\) is coercive, then the following also hold:
		\begin{enumerate}[resume]
		\item
			\(\seq{z^k}\) is bounded (in fact, this holds surely for arbitrary sampling criteria).
		\end{enumerate}
	\end{cor}
	
	\begin{cor}[linear convergence of \Cref{alg:Finito} under strong convexity]%
		Additionally to the assumptions of \Cref{thm:Finito:convergence}, suppose that \(g\) is convex and that each \(f_i\) is \(\mu_{f_i}\)-strongly convex.
		The following hold for the iterates generated by \Cref{alg:Finito}:
		\begin{itemize}[leftmargin=*,label={},itemindent=-0.5cm,labelsep=0pt,partopsep=0pt,parsep=0pt,listparindent=0pt,topsep=0pt]
		\item{\sc Randomized sampling:}
			under \Cref{ass:random},
			\begin{align*}
				\E[]{\varphi(z^k)-\min\varphi}
			{}\leq{} &
				(\varphi(x^{\rm init})-\min\varphi)
				(1-c)^k
			\\
				\tfrac12\E[]{\|z^k-x^\star\|^2}
			{}\leq{} &
				\frac{
					N(\varphi(x^{\rm init})-\min\varphi)
				}{
					\sum_i\mu_{f_i}
				}
				(1-c)^k
			\end{align*}
			holds for all \(k\in\N\), where \(c\) is as in \eqref{eq:cwc} and \(x^\star\coloneqq\argmin\varphi\).
			If the stepsizes \(\gamma_i\) and the sampling probabilities \(p_i\) are set as in \Cref{thm:random:linear}, then the tighter constant \(c\) as in \eqref{eq:cbc} is obtained.
		\item{\sc Shuffled cyclic or cyclic sampling:}
			under either sampling strategy \eqref{eq:ShufCyclicRule} or \eqref{eq:cyclicRule},
			\begin{align*}
				\varphi(z^{\nu N})-\min\varphi
			{}\leq{} &
				(\varphi(x^{\rm init})-\min\varphi)
				(1-c)^\nu
			\\
				\tfrac12\E[]{\|z^{\nu N}-x^\star\|^2}
			{}\leq{} &
				\frac{
					N(\varphi(x^{\rm init})-\min\varphi)
				}{
					\sum_i\mu_{f_i}
				}
				(1-c)^\nu
			\end{align*}
			holds surely for all \(\nu\in\N\), where \(c\) is as in \eqref{eq:linearShuffledCyclic}.
		\end{itemize}
	\end{cor}
	
	The next result follows from \Cref{thm:cyclic:global} once the needed properties of \(\Phi\) as in the umbrella formulation \eqref{eq:P} are shown to hold.
	
	\begin{cor}[global convergence of \Cref{alg:Finito}]\label{thm:Finito:global}%
		In the finite sum problem \eqref{eq:FSP}, suppose that \(\varphi\) has the KL property with exponent \(\theta\in(0,1)\) (as is the case when \(f_i\) and \(g\) are semialgebraic) and coercive, \(g\) is proper convex and lsc, and each \(f_i\) is \(L_{f_i}\)-Lipschitz differentiable, \(i\in[N]\).
		Then the sequence \(\seq{z^k}\) generated by \Cref{alg:Finito} with any essentially cyclic sampling strategy as in \Cref{ass:cyclic} converges surely to a stationary point for \(\varphi\).
		Moreover, if \(\theta\leq\nicefrac12\) then it converges at \(R\)-linear rate.
		\begin{proof}
			Function \(\Phi=F+G\) be as in \eqref{eq:FINITOG} clearly is coercive and satisfies \Cref{ass:basic}.
			In order to invoke \Cref{thm:cyclic:global} is suffices to show that there exists a constant \(c>0\) such that
			\begin{equation}\label{eq:Dist}
				\dist(0,\partial\Phi(\bm x))
			{}\geq{}
				c\dist(0,\partial\varphi(x))
			\quad
				\text{for all \(x\in\R^n\) and \(\bm x=(x,\dots,x)\),}
			\end{equation}
			as this will ensure that \(\Phi\) enjoys the KL property at \(\bm x^\star=(x^\star,\dots,x^\star)\) with same desingularizing function (up to a positive scaling).
			Notice that for \(x\in\R^n\) and \(\bm x=(x,\dots,x)\), one has
			\(
				\bm v\in\partial G(\bm x)
			\)
			iff
			\(
				\frac1N\sum_{i=1}^Nv_i
			{}\in{}
				\partial g(x)
			\).
			Since
			\(
				\partial\Phi(\bm x)
			{}={}
				\tfrac1N\mathop\times_{i=1}^N\nabla f_i(x_i)+\partial G(\bm x)
			\)
			and
			\(
				\partial\varphi(x)
			{}={}
				\tfrac1N\sum_{i=1}^N\nabla f_i(x)
				{}+{}
				\partial g(x)
			\),
			see \cite[Ex. 8.8(c) and Prop. 10.5]{rockafellar2011variational}, for \(x\in\R^n\) and denoting \(\bm x=(x,\dots,x)\) we have
			\begin{align*}
				\dist(0,\partial\varphi(x))
			{}\leq{} &
				\inf_{\bm v\in\partial G(\bm x)}{
					\left\|\textstyle
						\tfrac1N\sum_{i=1}^N\nabla f_i(x)
						{}+{}
						\tfrac1N\sum_{i=1}^Nv_i
					\right\|
				}
			\\
			{}\leq{} &
				\tfrac1N\inf_{\bm v\in\partial G(\bm x)}{
					\textstyle
					\sum_{i=1}^N\|\nabla f_i(x)+v_i\|
				}
			{}={}
				\tfrac1N\inf_{\bm u\in\partial\Phi(\bm x)}{
					\newnorm{\bm u}
				},
			\end{align*}
			where \(\newnorm{{}\cdot{}}\) is the norm in \(\R^{Nn}\) given by
			\(
				\newnorm{\bm w}=\sum_{i=1}^N\|w_i\|
			\).
			Inequality \eqref{eq:Dist} then follows by observing that
			\(
				\inf_{\bm u\in\partial\Phi(\bm x)}{
					\newnorm{\bm u}
				}
			\)
			is the distance of \(0\) from \(\partial\Phi(\bm x)\) in the norm \(\newnorm{{}\cdot{}}\), hence that \(\newnorm{{}\cdot{}}\leq c'\|{}\cdot{}\|\) for some \(c'>0\).
		\end{proof}
	\end{cor}

	\section{Nonconvex sharing problem}\label{sec:Sharing}

	In this section we consider the sharing problem \eqref{eq:SP}. 
	As discussed in \Cref{sec:Introduction}, \eqref{eq:SP} fits into the problem framework \eqref{eq:P} by simply letting \(G\coloneqq g \circ A\), where
	\(A\coloneqq[\I_n~\dots~\I_n]\in\R^{n\times nN}\). 
	By arguing as in \cite[Th. 6.15]{beck2017first} it can be shown that, when \(A\) has full row rank, the proximal mapping of $G=g\circ A$ is given by
	\begin{equation}\label{eq:sharingprox}
		\prox_G^{\Gamma^{-1}}(\bm u)
	{}={}
		\bm u+\Gamma\trans A(A\Gamma\trans A\,)^{-1}\left(\prox^{(A\Gamma\trans A\,)^{-1}}_g\left(A\bm u\right)-A\bm u\right).
	\end{equation}
	Since \(A\Gamma\trans A=\sum_{i=1}^N\gamma_i\) for the sharing problem \eqref{eq:SP},
	\begin{align*}
		\bm v
	{}\in{}
		\prox_G^{\Gamma^{-1}}(\bm u)
	~~\Leftrightarrow~~ &
		\bm v
	{}={}
		(u_1+\gamma_1w,\dots,u_N+\gamma_Nw)
	\\
	&
	\textstyle
		w
	{}\in{}
		\tilde\gamma^{-1}\left(\prox_{\tilde{\gamma}g}(\tilde u)-\tilde u\right),
	~~
		\tilde\gamma\coloneqq\sum_{i=1}^N\gamma_i,
	~~
		\tilde u\coloneqq \sum_{i=1}^Nu_i.
	\end{align*}
	Consequently general BC \Cref{alg:BC} when applied to the sharing problem \eqref{eq:SP} reduces to \Cref{alg:Sharing}. 
	
	\begin{algorithm} 
		\caption{Block-coordinate method for nonconvex sharing problem \eqref{eq:SP}}%
		\label{alg:Sharing}%
	\begin{algorithmic}[1]
	\item[{\sc Require}]
		\(
			 x_i^{\rm init}\in\R^{n}
		\),~
		\(\gamma_i\in(0,\nicefrac{N}{L_{f_i}})\),
		{\small \(i\in [N]\)}
	\Statex
		\(
			\tilde\gamma
		{}\coloneqq{}
			 \sum_{i=1}^N\gamma_i
		\),~~
		\(
			s_i
		{}={}
			 x_i^{\rm init}-\frac{\gamma_i}{N}\nabla f_i(x_i^{\rm init})
		\)~
		\(i\in [N]\),~~\(
			\tilde s
		{}={}
			\sum_{i=1}^N s_i
		\)
	\item[{\sc Repeat} until convergence]
		\State
			select a set of indices \(I\subseteq[N]\)
		\State $w \gets \tilde{\gamma}^{-1}(\prox_{\tilde{\gamma}g}(\tilde s)-\tilde s)$
	
		\For{ \(i\in I\) }
		\State 
		\(
		v_i {}\gets{} s_i + \gamma_i w  - \tfrac{\gamma_i}{N} \nabla f_i(s_i + \gamma_i w )
		\)
		\State
			update~~
			\(
				\tilde s
			{}\gets{}
				\tilde s+(v_i-s_i)
			\)
			~~and~~
			\(
				s_i \gets v_i
			\)
			\EndFor
			\item[{\sc Return}]
				$\bm z=(s_1 + \gamma_1 w ,\dots,s_N + \gamma_N w)$ with $w\in\tilde{\gamma}^{-1}(\prox_{\tilde{\gamma}g}(\tilde s)-\tilde s)$
	\end{algorithmic}
	
	
	\end{algorithm}
	
	\begin{rem}[generalized sharing constraint]
		Another notable instance of $G=g\circ A$ well suited for the BC framework of \Cref{alg:BC} is when \(g=\indicator_{\set0}\) and \(A=[A_1~\dots~A_N]\), \(A_i\in\R^{n\times n_i}\) such that $A$ is full rank.
		This models the generalized sharing problem
		\[
			\minimize_{\bm x\in\R^{\sum_in_i}}{\textstyle
				\tfrac1N\sum_{i=1}^Nf_i(x_i)
			}
		\quad\stt{}\textstyle
			\sum_{i=1}^NA_ix_i=0.
		\]
		In this case \eqref{eq:sharingprox} simplifies to
		\[
			\left(\prox_G^{\Gamma^{-1}}(\bm u)\right)_i
		{}={}
			u_i-\gamma_i\trans{A_i}\mathcal A^{-1}\sum_{i=1}^NA_iu_i, 
		\]
		where  $\mathcal A\coloneqq A\Gamma\trans A$ can be factored offline and \(\sum_{i=1}^NA_ix_i\) can be updated in an incremental fashion in the same spirit of \Cref{alg:Sharing}.
	\end{rem}
	
	The convergence results for \Cref{alg:Sharing} summarized below fall as special cases of those in  \Cref{sec:convergence}.
	
	\begin{cor}[convergence of \Cref{alg:Sharing}]\label{thm:sharing:convergence}%
		In the sharing problem \eqref{eq:SP}, suppose that \(\argmin\Phi\) is nonempty, \(g\) is proper and lsc, and each \(f_i\) is \(L_{f_i}\)-Lipschitz differentiable, \(i\in[N]\). Consider the sequences $\seq{w^k}$ and $\seq{\bm s^k}$ generated by \Cref{alg:Sharing} and let $\seq{\bm z^k}=\seq{s_1^k + \gamma_1 w^k ,\dots,s_N^k + \gamma_N w^k}$.  
		Then, the following hold almost surely (resp. surely) with randomized sampling strategy as in \Cref{ass:random} (resp. with any essentially cyclic sampling strategy and $g$ convex as required in \Cref{ass:cyclic}):
		\begin{enumerate}
		\item
			the sequence \(\seq{\Phi(\bm z^k)}\) converges to a finite value \(\Phi_\star\leq\Phi(\bm x^{\rm init})\);
		\item
			all cluster points of the sequence \(\seq{\bm z^k}\) are stationary and on which \(\Phi\) equals \(\Phi_\star\).
		\end{enumerate}
		If, additionally, \(\Phi\) is coercive, then the following also hold:
		\begin{enumerate}[resume]
		\item
			\(\seq{\bm z^k}\) is bounded (in fact, this holds surely for arbitrary sampling criteria).
		\end{enumerate}
	\end{cor}
	
	\begin{cor}[linear convergence of \Cref{alg:Sharing} under strong convexity]\label{cor:RLinSharing}%
		Additionally to the assumptions of \Cref{thm:sharing:convergence}, suppose that \(g\) is convex and that each \(f_i\) is \(\mu_{f_i}\)-strongly convex.
		The following hold:
		\begin{itemize}[leftmargin=*,label={},itemindent=-0.5cm,labelsep=0pt,partopsep=0pt,parsep=0pt,listparindent=0pt,topsep=0pt]
		\item{\sc Randomized sampling:}
			under \Cref{ass:random},
			\begin{align*}
				\E[]{\Phi(\bm z^k)-\min\Phi}
			{}\leq{} &
				\bigl(\Phi(\bm x^{\rm init})-\min\Phi\bigr)(1-c)^k
			\\
				\tfrac12\E[]{\|\bm z^k-\bm x^\star\|^2_{\mu_F}}
			{}\leq{} &
				\bigl(\Phi(\bm x^{\rm init})-\min\Phi\bigr)(1-c)^k
			\end{align*}
			holds for all \(k\in\N\), where \(\bm x^\star\coloneqq\argmin\Phi\),
			\(
				\mu_F
			{}\coloneqq{} \tfrac1N
				\blockdiag\bigl(\mu_{f_1}\I_{n_1},\dots\mu_{f_n}\I_{n_N}\bigr)
			\),
			and \(c\) is as in \eqref{eq:cwc}.
			If the stepsizes \(\gamma_i\) and the sampling probabilities \(p_i\) are set as in \Cref{thm:random:linear}, then the tighter constant \(c\) as in \eqref{eq:cbc} is obtained.
		\item{\sc Shuffled cyclic or cyclic sampling:}
			under either sampling strategy \eqref{eq:ShufCyclicRule} or \eqref{eq:cyclicRule},
			\begin{align*}
				\Phi(\bm z^{N\nu})-\min\Phi
			{}\leq{} &
				\bigl(\Phi(\bm x^{\rm init})-\min\Phi\bigr)(1-c)^\nu
			\\
				\tfrac12\|\bm z^{N\nu}-\bm x^\star\|^2_{\mu_F}
			{}\leq{} &
				\bigl(\Phi(\bm x^{\rm init})-\min\Phi\bigr)(1-c)^\nu
			\end{align*}
			holds surely for all \(\nu\in\N\), where \(c\) is as in \eqref{eq:linearShuffledCyclic}.
		\end{itemize}
	\end{cor}
	
	We conclude with an immediate consequence of \Cref{thm:cyclic:global} that shows that (strong) convexity is in fact not necessary for global or linear convergence to hold.
	
	\begin{cor}[global and linear convergence of \Cref{alg:Sharing}]\label{thm:Sharing:global}%
		In problem \eqref{eq:SP}, suppose that \(\Phi\) has the KL property with exponent \(\theta\in(0,1)\) (as is the case when \(g\) and \(f_i\) are semialgebraic) and is coercive, \(g\) is proper convex lsc, and each \(f_i\) is \(L_{f_i}\)-Lipschitz differentiable, \(i\in[N]\).
		Then the sequence $\seq{\bm z^k}$ as defined in \Cref{thm:sharing:convergence} with any essentially cyclic sampling strategy as in \Cref{ass:cyclic} converges surely to a stationary point for \(\Phi\).
		Moreover, if \(\theta\leq\nicefrac12\) it converges with \(R\)-linear rate.
	\end{cor}

\ifaccel
	\section{Accelerated block-coordinate proximal gradient}
	The work \cite{allen2016even} introduced a coordinate descent method for smooth convex minimization, in which each coordinate is randomly sampled according to an ad hoc probability distribution that provably leads to a remarkable speed up with respect to uniform sampling strategies.
	The unified analysis of BC-algorithms and the analytical tool introduced in this paper, the forward backward envelope function, allow the extention of this approach to nonsmooth convex minimization of the form \eqref{eq:P},
	where functions \(f_i\) are convex quadratic and \(G\) is convex but possibly nonsmooth:
	\begin{ass}[requirements for the fast BC-\Cref{alg:Fast}]\label{ass:Fast}%
		In problem \eqref{eq:P}, \(\func{G}{\R^{\sum_in_i}}{\Rinf}\) is proper convex and lsc, and
		\(f_i(x_i)\coloneqq\tfrac12\trans{x_i}H_ix_i+\trans{q_i}x_i\) is convex quadratic, with \(L_{f_i}\coloneqq\lambda_{\rm max}(H_i)\) and \(\mu_{f_i}\coloneqq\lambda_{\rm min}(H_i)\geq0\),
		\(i\in[N]\).
	\end{ass}
	
	Let $U_i\in \R^{{\sum_{i=1}^N n_i}\times n_i}$ denote the $i$-th block column of the identity matrix so that for a vector $v\in \R^{n_i}$
	\begin{equation}\label{eq:U}
		U_iv= (0,\dots, 0, \!\overbracket{\,v\,}^{\mathclap{i\text{-th}}}\!, 0, \dots, 0).
	\end{equation}
	
	The accelerated BC scheme based on \cite{allen2016even} (for both strongly convex and convex cases) is given in \Cref{alg:Fast}.
	Similarly to the approach of \cite{patrinos2014douglas} where an accelerated Douglas-Rachford algorithm is proposed, in order to derive \Cref{alg:Fast} we consider the scaled problem \( \minimize_{\tilde{\bm x}} \FBEC(\tilde{\bm x})\) where $\FBEC \coloneqq \FBE\circ Q^{-1/2}$, and $Q$ is the symmetric positive definite matrix 
	\begin{equation} \label{eq:QQ}
		Q
	{}\coloneqq{}
		\blockdiag(Q_1,\dots,Q_N)\succ0
	\quad \text{with } Q_i
	{}\coloneqq{}
		\gamma_i^{-1}\I-\tfrac{1}{N}H_i\in\R^{n_i\times n_i},~i\in[N]. 
	\end{equation}
	As detailed in \Cref{thm:convex}, whenever \Cref{ass:Fast} is satisfied $\FBEC$ is a convex Lipschitz-differentiable function, and its gradient is given by $\nabla \FBEC(\tilde{\bm x}) = Q^{1/2}(\bm x-\prox_G^{\Gamma^{-1}}(\bm x-\Gamma\nabla F(\bm x)))$ where $\bm x=Q^{-1/2}\tilde{\bm x}$.
	Note that, based on \Cref{thm:convex}, \(\FBEC\) is \(1\)-smooth along the \(i\)-th block (in the notation of \cite{allen2016even}, \(L_i=1\), $S_\alpha=N$, and \(p_i=\nicefrac1N\)).
	Hence the parameters of the algorithm simplify substantially resulting in uniform sampling. 
	Moreover, when functions \(f_i\) are \(\mu_{f_i}\)-strongly convex, by \Cref{thm:convex} \(\FBEC\) is $\sigma$-strongly convex with $\sigma= \frac{1}{N} \min_{i\in [N]} \{\gamma_i\mu_{f_i}\}$. 
	
	\Cref{alg:Fast} is obtained by applying the fast BC to this problem and scaling the variables by $Q^{-1/2}$.
	Specifically, the update rule as in \cite{allen2016even} reads
	\[
		\begin{cases}[ r @{{}={}} l ]
			\tilde{\bm x}^+ & \tau\tilde{\bm w}+(1-\tau)\tilde{\bm y}
		\\
			\tilde{\bm y}^+ & \tilde{\bm x}-U_i\trans{U_i}\nabla\FBEC(\tilde{\bm x}^+)
			{}={}
			\tilde{\bm x}-U_iQ_i^{\nicefrac12}(x_i^+-z_i^+)
		\\
			\tilde{\bm w}^+ & \tfrac{1}{1+\eta\sigma}(\tilde{\bm w}+\eta\sigma\tilde{\bm x}^+-N\eta U_i\trans{U_i}\nabla\FBEC(\tilde{\bm x}^+))
			{}={}
			\tfrac{1}{1+\eta\sigma}(\tilde{\bm w}+\eta\sigma\tilde{\bm x}^+-N\eta U_iQ_i^{\nicefrac12}(x_i^+-z_i^+)),
		\end{cases}
	\]
	where \(\bm z^+=\prox_G^{\Gamma^{-1}}(\bm x^+-\Gamma\nabla F(\bm x^+))\).
	Since \(Q^{-\nicefrac12}U_iQ_i^{\nicefrac12}=U_i\), premultiplying by \(Q^{-\nicefrac12}\) yields
	\[
		\begin{cases}[ r @{{}={}} l ]
			\bm x^+ & \tau\bm z+(1-\tau)\bm y
		\\
			\bm z^+ & \prox_G^{\Gamma^{-1}}(\bm x^+-\Gamma\nabla F(\bm x^+))
		\\
			\bm y^+ & \bm x+U_i(z_i^+-x_i^+)
		\\
			\bm w^+ & \tfrac{1}{1+\eta\sigma}(\bm w+\eta\sigma\bm x^++N\eta U_i(z_i^+-x_i^+)).
		\end{cases}
	\]
	For computational efficiency, vectors $\Gamma\nabla F(\bm x^k)$ and $\Gamma\nabla F(\bm w^k)$ are stored in variables $\bm r^k$ and $\bm v^k$ and updated recursively using the fact that gradients are affine, in such a way that each iteration requires only the evaluation of the sampled gradient (see \Cref{state:d}).
	For similar reasons, in \Cref{alg:Fast} the iterates start with the $\bm y$-update rather than the $\bm x$-update as in \cite{allen2016even}.     
	Moreover, in the same spirit of \Cref{alg:BC} this accelerated variant can be implemented efficiently whenever the individual blocks of \(\bm z^+\) can be computed efficiently, similarly to the cases discussed in \Cref{sec:Finito,sec:Sharing}.
	
	\begin{algorithm}
		\caption{Accelerated block-coordinate proximal gradient for problem \eqref{eq:P} under \Cref{ass:Fast}}%
		\label{alg:Fast}%
	\begin{algorithmic}[1]
	\item[{\sc Require}]
		\(\bm x^0\in\R^{\sum_in_i}\),~
		\(\gamma_i\in(0,\nicefrac{N}{L_{f_i}})\),~$i\in[N]$,
		\(
			\sigma
		{}\coloneqq{}
			\frac{1}{N} \min_{i\in [N]} \{\gamma_i\mu_{f_i}\}
		\)
	\State\label{state:sigma_beta}%
		{\bf if} \(\sigma=0\),~~{\bf then}~~%
		\(\eta = \nicefrac{1}{N^2}\)
		~~{\bf otherwise}~~%
		set
		\(
			\tau
		{}={}
			\frac{2}{1+\sqrt{1+\nicefrac{4N^2}{\sigma}}}
		\),
		\(
			\eta
		{}={}
			\frac{1}{\tau N^2}
		\)~~%
			{\bf end if}
	\State
		\(
			\bm w^0
		{}={}
			\bm x^0
		\),~
		\(
			(\bm v^0,\bm r^0)
		{}={}
			(\Gamma\nabla F(\bm x^0),\Gamma\nabla F(\bm x^0))
		\),~
		\(
			\bm z^{0}
			{}={}
			\prox_G^{\Gamma^{-1}}\bigl(\bm x^{0}-\bm r^{0}\bigr)
		\)
		\def\myVar#1{%
			\fillwidthof[l]{\bm x^{k+1}}{#1}%
		}%
	\For{ \(k=0,1,\dots\) }
		\State
			sample \(i\in[N]\) uniformly
		\State\label{state:d}%
			\(
				\myVar{\bm y^{k+1}}
			{}\gets{}
				\bm x^{k}+U_{i}\bigl(z_{i}^{k}-x_{i}^{k}\bigr)
			\),\quad%
			\(
				d
			{}\gets{}
				 \tfrac{\gamma_i}{N}\nabla f_i(z_i^{k})-r_i^{k} 
			\)
			\State 
			\(
				\myVar{\bm v^{k+1}}
			{}={}
				\frac{1}{1+\eta\sigma}\Bigl(\bm v^{k}+\eta\sigma\bm r^{k}+N{\eta} U_id\Bigr)
			\),\quad
			\(
				\bm w^{k+1}
			{}={}
				\frac{1}{1+\eta\sigma}\Bigl(\bm w^{k}+\eta\sigma\bm x^{k}+{N\eta}U_{i}\bigl(z_{i}^{k}-x_{i}^{k}\bigr)\Bigr)
			\)
		\State{\bf if}~~\(\sigma=0\),~~{\bf then}~~%
			\(
				\eta
			{}\gets{}
				\frac{k+3}{2N^2}
			\),~~%
			\(
				\tau
			{}\gets{}
				\frac{2}{k+3}
			\);%
			~~{\bf end if}%
		\State\label{state:FBEgrad}%
			\(
				\myVar{\bm x^{k+1}}
			{}={}
				\tau \bm w^{k+1}+(1-\tau)\bm y^{k+1}
			\),\quad
			\(
				\bm r^{k+1}
			{}={}
				\tau \bm v^{k+1}+(1-\tau)(\bm r^{k}+U_i d)
			\)%
			\State
			\(
				\myVar{\bm z^{k+1}}
			{}={}
				\prox_G^{\Gamma^{-1}}\bigl(\bm x^{k+1}-\bm r^{k+1}\bigr)
			\)
	\EndFor{}
	\end{algorithmic}

	\end{algorithm}
	
	The convergence rate results follow directly from those of \cite{allen2016even} with parameters \(L_i=1\) and $S_\alpha=N$ as described above.
	
	\begin{thm}[convergence rates of \Cref{alg:Fast}]
		Suppose that \Cref{ass:basic,ass:Fast} are satisfied.
		Then, the iterates generated by \Cref{alg:Fast} satisfy
		\[
			\mathbb{E}\left[{\FBE(\bm y^k)-\min \Phi}\right]
		{}\leq{}
			\frac{2N^2\|\bm x^0 - \bm x^\star\|^2_Q}{(k+1)^2},
		\]
		where $Q$ is as in \eqref{eq:QQ}. 
		Moreover, in the strongly convex case ($\sigma= \frac{1}{N} \min_{i\in [N]} \{\gamma_i\mu_{f_i}\}>0$)
		\[
			\mathbb{E}\left[{\FBE(\bm y^k)-\min  \Phi}\right]
		{}\leq{}
			O(1)
			(1-c)^k\left(
				\FBE(\bm x^0)-\min  \Phi
			\right)
		\quad\text{where}\quad
		\textstyle
			c
		{}={}
			\left(
				\frac12
				{}+{}
				\sqrt{
					\frac14
					{}+{}
					\frac{N^2}{\sigma}
				}
			\right)^{-1}.
		\]
	\end{thm}
	
	Note that in the strongly convex case it follows from \Cref{thm:strconcost} that the distance from the solution decreases \(R\)-linearly as
	\[
		\E[]{ \|\bm y^{k}-\bm x^\star\|^2_M} 
	{}\leq{}
		O(1)\left(1-c\right)^k
		\left(\FBE(\bm x^0)-\min\Phi\right),
	\]
	where $M$ is as in \Cref{thm:strconcost}.
\fi

	\section{Conclusions}\label{sec:Conclusions}

	We presented a general block-coordinate forward-backward algorithm for minimizing the sum of a separable smooth and a nonseparable nonsmooth function, both allowed to be nonconvex.
	The framework is general enough to encompass regularized finite sum minimization and sharing problems, and leads to (a generalization of) the Finito/MISO algorithm \cite{defazio2014finito,mairal2015incremental} with new convergence results and with another novel incremental-type algorithm.
	The forward-backward envelope is shown to be a particularly suitable Lyapunov function for establishing convergence: additionally to enjoying favorable continuity properties, \emph{sure} descent (as opposed to in expectation) occurs along the iterates.
	Possible future developments include extending the framework to account for a nonseparable smooth term, for instance by ``quantifying the strength of coupling'' between blocks of variables as in \cite[\S7.5]{bertsekas1989parallel}.

\ifarxiv
	\clearpage
\fi

	\begin{appendix}

		\section{The key tool: the forward-backward envelope}\label{sec:appendix}
	This appendix contains some proofs and auxiliary results omitted in the main body.
	We begin by observing that, since \(F\) and \(-F\) are 1-smooth in the metric induced by
	\(
		\Lambda_F\coloneqq\tfrac1N\blockdiag(L_{f_1}\I_{n_1},\dots,L_{f_N}\I_{n_N})
	\),
	one has
	\begin{equation}\label{eq:Lip}
		F(\bm x)+\innprod{\nabla F(\bm x)}{\bm w-\bm x}
		{}-{}
		\tfrac12\|\bm w-\bm x\|_{\Lambda_F}^2
	{}\leq{}
		F(\bm w)
	{}\leq{}
		F(\bm x)+\innprod{\nabla F(\bm x)}{\bm w-\bm x}
		{}+{}
		\tfrac12\|\bm w-\bm x\|_{\Lambda_F}^2
	\end{equation}
	for all \(\bm x,\bm w\in\R^{\sum_in_i}\), see \cite[Prop. A.24]{bertsekas2016nonlinear}.
	Let us denote
	\[
		\M(\bm w,\bm x)
	{}\coloneqq{}
		F(\bm x)+\innprod{\nabla F(\bm x)}{\bm w-\bm x}
		{}+{}
		G(\bm w)
		{}+{}
		\tfrac12\|\bm w-\bm x\|_{\Gamma^{-1}}^2
	\]
	the quantity being minimized (with respect to \(\bm w\)) in the definition \eqref{eq:FBE} of the FBE.
	It follows from \eqref{eq:Lip} that
	\begin{equation}\label{eq:bounds}
		\Phi(\bm w)
		{}+{}
		\tfrac12\|\bm w-\bm x\|^2_{\Gamma^{-1}-\Lambda_F}
	{}\leq{}
		\M(\bm w,\bm x)
	{}\leq{}
		\Phi(\bm w)
		{}+{}
		\tfrac12\|\bm w-\bm x\|^2_{\Gamma^{-1}+\Lambda_F}
	\end{equation}
	holds for all \(\bm x,\bm w\in\R^{\sum_in_i}\).
	In particular, \(\M\) is a \emph{majorizing model} for \(\Phi\), in the sense that \(\M(\bm x,\bm x)=\Phi(\bm x)\) and \(\M(\bm w,\bm x)\geq\Phi(\bm w)\) for all \(\bm x,\bm w\in\R^{\sum_in_i}\).
	In fact, as explained in \Cref{sec:FBE}, while a \(\Gamma\)-forward-backward step \(\bm z\in\T(\bm x)\) amounts to evaluating a minimizer of \(\M({}\cdot{},\bm x)\), the FBE is defined instead as the minimization value, namely \(\FBE(\bm x)=\M(\bm z,\bm x)\) where \(\bm z\) is any element of \(\T(\bm x)\).

			\subsection{Proofs of \texorpdfstring{\Cref{sec:FBE}}{\S\ref*{sec:FBE}}}\label{sec:proofs:FBE}
	\begin{appendixproof}{thm:osc}
		For \(\bm x^\star\in\argmin\Phi\) it follows from \eqref{eq:Lip} that
		\[
			\min\Phi
		{}\leq{}
			F(\bm x)
			{}+{}
			G(\bm x)
		{}\leq{}
			G(\bm x)
			{}+{}
			F(\bm x^\star)
			{}+{}
			\innprod{\nabla F(\bm x^\star)}{\bm x-\bm x^\star}
			{}+{}
			\tfrac12\|\bm x^\star-\bm x\|_{\Lambda_F}^2.
		\]
		Therefore, \(G\) is lower bounded by a quadratic function with quadratic term \(-\tfrac12\|{}\cdot{}\|_{\Lambda_F}^2\), and thus is prox-bounded in the sense of \cite[Def. 1.23]{rockafellar2011variational}.
		The claim then follows from \cite[Th. 1.25 and Ex. 5.23(b)]{rockafellar2011variational} and the continuity of the forward mapping \(\Fw{}\).
	\end{appendixproof}

	\begin{appendixproof}{thm:FBEineq}%
		Local Lipschitz continuity
		\ifarxiv
			of the FBE
		\fi
		follows from \eqref{eq:FBEMoreau} in light of \Cref{thm:osc} and \cite[Ex. 10.32]{rockafellar2011variational}.
		\begin{proofitemize}
		\item\ref{thm:leq}~
			Follows by replacing \(\bm w=\bm x\) in \eqref{eq:FBE}.
		\item\ref{thm:geq}~
			Directly follows from \eqref{eq:bounds} and the identity \(\FBE(\bm x)=\M(\bm z,\bm x)\) for \(\bm z\in\T(\bm x)\).
		\item\ref{thm:strconcost}~
			By strong convexity, denoting \(\Phi_\star\coloneqq\min\Phi\), we have 	
			\[
				\Phi_\star
			{}\leq{}
				\Phi(\bm z)-\tfrac12\|\bm z-\bm x^\star\|_{\mu_F}^2
			{}\leq{}
				\FBE(\bm x)
				{}-{}
				\tfrac12\|\bm z-\bm x^\star\|_{\mu_F}^2
			\]
			where the second inequality follows from \Cref{thm:geq}.
		\qedhere
		\end{proofitemize}
	\end{appendixproof}
	
	\begin{appendixproof}{thm:FBEmin}
		\begin{proofitemize}
		\item\ref{thm:min} and \ref{thm:argmin}~
			It follows from \Cref{thm:leq} that \(\inf\FBE\leq\min\Phi\).
			Conversely, let \(\seq{\bm x^k}\) be such that \(\FBE(\bm x^k)\to\inf\FBE\) as \(k\to\infty\), and for each \(k\) let \(\bm z^k\in\T(\bm x^k)\).
			It then follows from \Cref{thm:leq,thm:geq} that
			\[
				\inf\FBE
			{}\leq{}
				\min\Phi
			{}\leq{}
				\liminf_{k\to\infty}\Phi(\bm z^k)
			{}\leq{}
				\liminf_{k\to\infty}\FBE(\bm x^k)
			{}={}
				\inf\FBE,
			\]
			hence \(\min\Phi=\inf\FBE\).
	%
			Suppose now that \(\bm x\in\argmin\Phi\) (which exists by \Cref{ass:basic}); then it follows from \Cref{thm:geq} that \(\T(\bm x)=\set{\bm x}\) (for otherwise another element would belong to a lower level set of \(\Phi\)).
			Combining with \Cref{thm:leq} with \(\bm z=\bm x\) we then have
			\[
				\min\Phi
			{}={}
				\Phi(\bm z)
			{}\leq{}
				\FBE(\bm x)
			{}\leq{}
				\Phi(\bm x)
			{}={}
				\min\Phi.
			\]
			Since \(\min\Phi=\inf\FBE\), we conclude that \(\bm x\in\argmin\FBE\), and that in particular \(\inf\FBE=\min\FBE\).
			Conversely, suppose \(\bm x\in\argmin\FBE\) and let \(\bm z\in\T(\bm x)\).
			By combining \Cref{thm:leq,thm:geq} we have that \(\bm z=\bm x\), that is, that \(\T(\bm x)=\set{\bm x}\).
			It then follows from \Cref{thm:geq} and assert \ref{thm:min} that
			\[
				\Phi(\bm x)
			{}={}
				\Phi(\bm z)
			{}\leq{}
				\FBE(\bm x)
			{}={}
				\min\FBE
			{}={}
				\min\Phi,
			\]
			hence \(\bm x\in\argmin\Phi\).
		\item\ref{thm:LB}~
			Due to \Cref{thm:leq}, if \(\FBE\) is level bounded clearly so is \(\Phi\).
			Conversely, suppose that \(\FBE\) is not level bounded.
			Then, there exist \(\alpha\in\R\) and \(\seq{\bm x^k}\subseteq\lev_{\leq\alpha}\FBE\) such that \(\|\bm x^k\|\to\infty\) as \(k\to\infty\).
			Let \(\lambda=\min_i\set{\gamma_i^{-1}-L_{f_i}N^{-1}}>0\), and for each \(k\in\N\) let \(\bm z^k\in\T(\bm x^k)\).
			It then follows from \Cref{thm:geq} that
			\[
				\min\Phi
			{}\leq{}
				\Phi(\bm z^k)
			{}\leq{}
				\FBE(\bm x^k)
				{}-{}
				\tfrac\lambda2\|\bm x^k-\bm z^k\|^2
			{}\leq{}
				\alpha
				{}-{}
				\tfrac\lambda2\|\bm x^k-\bm z^k\|^2,
			\]
			hence \(\seq{\bm z^k}\subseteq\lev_{\leq\alpha}\Phi\) and
			\(
				\|\bm x^k-\bm z^k\|^2
			{}\leq{}
				\tfrac2\lambda(\alpha-\min\Phi)
			\).
			Consequently, also the sequence \(\seq{\bm z^k}\subseteq\lev_{\leq\alpha}\Phi\) is unbounded, proving that \(\Phi\) is not level bounded.
		\qedhere
		\end{proofitemize}
	\end{appendixproof}

			\subsection{Further results}\label{sec:auxiliary}
	This section contains a list of auxiliary results invoked in the main proofs of \Cref{sec:convergence}.
	
	\begin{lem}\label{thm:critical}%
		Suppose that \Cref{ass:basic} holds, and let two sequences \(\seq{\bm u^k}\) and \(\seq{\bm v^k}\) satisfy \(\bm v^k\in\T(\bm u^k)\) for all \(k\) and be such that both converge to a point \(\bm u^\star\) as \(k\to\infty\).
		Then, \(\bm u^\star\in\T(\bm u^\star)\), and in particular \(0\in\hat\partial\Phi(\bm u^\star)\).
		\begin{proof}
			Since \(\nabla F\) is continuous, it holds that \(\Fw{\bm u^k}\to\Fw{\bm u^\star}\) as \(k\to\infty\).
			From outer semicontinuity of \(\prox_G^{\Gamma^{-1}}\) \cite[Ex. 5.23(b)]{rockafellar2011variational} it then follows that
			\[
				\bm u^\star
			{}={}
				\lim_{k\to\infty}
				\bm v^k
			{}\in{}
				\limsup_{k\to\infty}
				\prox_G^{\Gamma^{-1}}(\Fw{\bm u^k})
			{}\subseteq{}
				\prox_G^{\Gamma^{-1}}(\Fw{\bm u^\star})
			{}={}
				\T(\bm u^\star),
			\]
			where the limit superior is meant in the Painlevé-Kuratowski sense, cf. \cite[Def. 4.1]{rockafellar2011variational}.
			The optimality conditions defining \(\prox_G^{\Gamma^{-1}}\) \cite[Th. 10.1]{rockafellar2011variational} then read
			\begin{align*}
				0
			{}\in{} &
				\hat\partial\left(
					G+\tfrac12\|{}\cdot{}-(\Fw{\bm u^\star})\|_{\Gamma^{-1}}^2
				\right)(\bm u^\star)
			{}={}
				\hat\partial G(\bm u^\star)
				{}+{}
				\Gamma^{-1}\left(
					\bm u^\star - (\Fw{\bm u^\star})
				\right)
			\\
			{}={} &
				\hat\partial G(\bm u^\star)
				{}+{}
				\nabla F(\bm u^\star)
			{}={}
				\hat\partial\Phi(\bm u^\star),
			\end{align*}
			where the first and last equalities follow from \cite[Ex. 8.8(c)]{rockafellar2011variational}.
		\end{proof}
	\end{lem}
	\begin{lem}%
		Suppose that \Cref{ass:basic} holds and that function \(G\) is convex.
		Then, the following hold:
		\begin{enumerate}
		\item\label{thm:FNE}
			\(\prox_G^{\Gamma^{-1}}\) is (single-valued and) firmly nonexpansive (FNE) in the metric $\|{}\cdot{}\|_{\Gamma^{-1}}$; namely,
			\[
				\|
					\prox_G^{\Gamma^{-1}}(\bm u)
					{}-{}
					\prox_G^{\Gamma^{-1}}(\bm v)
				\|_{\Gamma^{-1}}^2
			{}\leq{}
				\innprod{
					\prox_G^{\Gamma^{-1}}(\bm u)
					{}-{}
					\prox_G^{\Gamma^{-1}}(\bm v)
				}{
					\Gamma^{-1}(\bm u-\bm v)
				}
			{}\leq{}
				\|
					\bm u
					{}-{}
					\bm v
				\|_{\Gamma^{-1}}^2
			\quad\forall\bm u,\bm v;
			\]
		\item\label{thm:MoreauGrad}
			the Moreau envelope \(G^{\Gamma^{-1}}\) is differentiable with \(\nabla G^{\Gamma^{-1}}=\Gamma^{-1}(\id-\prox_G^{\Gamma^{-1}})\);
		\item\label{thm:subdiffdist}
			for every \(\bm x\in\R^{\sum_in_i}\) it holds that
			\(
				\dist(0,\partial\FBE(\bm x))
			{}\leq{}
				\tfrac{
					N+\max_i\set{\gamma_iL_{f_i}}
				}{
					N\min_i\set{\sqrt{\gamma_i}}
				}
				\|\bm x-\T(\bm x)\|_{\Gamma^{-1}}
			\);
		\item\label{thm:TLip}\label{thm:contractive}
			\(\T\) is \(L_{\bf T}\)-Lipschitz continuous in the metric $\|{}\cdot{}\|_{\Gamma^{-1}}$ for some \(L_{\bf T}\geq0\); if in addition \(f_i\) is \(\mu_{f_i}\)-strongly convex, \(i\in[N]\), then \(L_{\bf T}\leq 1-\delta\) for \(\delta=\frac1N\min_{i\in[N]}\set{\gamma_i\mu_{f_i}}\).
		\end{enumerate}
		\begin{proof}
			\begin{proofitemize}
			\item\ref{thm:FNE} and \ref{thm:MoreauGrad}~
				See \cite[Prop.s 12.28 and 12.30]{bauschke2017convex}.
			\item\ref{thm:subdiffdist}
				Let \(D\subseteq\R^{\sum_in_i}\) be the set of points at which \(\nabla F\) is differentiable.
				From the chain rule of differentiation applied to the expression \eqref{eq:FBEMoreau} and using assert \ref{thm:MoreauGrad}, we have that \(\FBE\) is differentiable on \(D\) with gradient
				\[
					\nabla\FBE(\bm x)
				{}={}
					\bigl[
						\I-\Gamma\nabla^2F(\bm x)
					\bigr]
					\Gamma^{-1}
					\bigl[
						\bm x-\T(\bm x)
					\bigr]
					\quad
					\forall\bm x\in D.
				\]
				Since \(D\) is dense in \(\R^{\sum_in_i}\) owing to Lipschitz continuity of \(\nabla F\), we may invoke \cite[Th. 9.61]{rockafellar2011variational} to infer that \(\partial\FBE(\bm x)\) is nonempty for every \(\bm x\in\R^{\sum_in_i}\) and
				\[
					\partial\FBE(\bm x)
				{}\supseteq{}
					\partial_B\FBE(\bm x)
				{}={}
					\bigl[
						\I-\Gamma\partial_B\nabla F(\bm x)
					\bigr]
					\Gamma^{-1}
					\bigl[
						\bm x-\T(\bm x)
					\bigr]
				{}={}
					\bigl[
						\Gamma^{-1}-\partial_B\nabla F(\bm x)
					\bigr]
					\bigl[
						\bm x-\T(\bm x)
					\bigr],
				\]
				where \(\partial_B\) denotes the (set-valued) Bouligand differential \cite[\S7.1]{facchinei2003finite}.
				The claim now follows by observing that
				\(
					\partial_B\nabla F(\bm x)
				{}={}
					\tfrac1N\blockdiag(\partial_B\nabla f_1(x_1),\dots,\partial_B\nabla f_N(x_N))
				\)
				and that each element of \(\partial_B\nabla f_i(x_i)\) has norm bounded by \(L_{f_i}\).
			\item\ref{thm:TLip}~
				Lipschitz continuity follows from assert \ref{thm:FNE} together with the fact that Lipschitz continuity is preserved by composition.
				Suppose now that \(f_i\) is \(\mu_{f_i}\)-strongly convex, \(i\in[N]\).
				By \cite[Thm 2.1.12]{nesterov2013introductory} for all $x_i,y_i\in\R^{n_i}$ 
				\begin{equation}\label{eq:smoothStrcvx}
					\langle\nabla f_i(x_i)-\nabla f_i(y_i),x_i-y_i\rangle\geq\tfrac{\mu_{f_i}L_{f_i}}{\mu_{f_i}+L_{f_i}}\|x_i-y_i\|^2+\tfrac1{\mu_{f_i}+L_{f_i}}\|\nabla f_i(x_i)-\nabla f_i(y_i)\|^2. 
				\end{equation}
				For the forward operator we have 
				\begin{align*}
				&
					\|
						(\id-\tfrac{\gamma_i}{N}\nabla f_i)(x_i)
						{}-{}
						(\id-\tfrac{\gamma_i}{N}\nabla f_i)(y_i)
					\|^2
				\\
				{}={} &
					\|x_i-y_i\|^2
					{}+{}
					\tfrac{\gamma_i^2}{N^2}
					\|\nabla f_i(x_i)-\nabla f_i(y_i)\|^2
					{}-{}
					\tfrac{2\gamma_i}{N}
					\innprod{x_i-y_i}{\nabla f_i(x_i)-\nabla f_i(y_i)}
				\\
				\overrel[\leq]{\eqref{eq:smoothStrcvx}}{} &
					\Bigl(
						1-\tfrac{\gamma_i^2\mu_{f_i}L_{f_i}}{N^2}
					\Bigr)
					\|x_i-y_i\|^2
					{}-{}
					\tfrac{\gamma_i}{N}
					\Bigl(
						2-\tfrac{\gamma_i}{N}(\mu_{f_i}+L_{f_i})
					\Bigr)
					\innprod{\nabla f_i(x_i)-\nabla f_i(y_i)}{x_i-y_i}
				\\
				{}\leq{} &
					\left(1-\tfrac{\gamma_i^2\mu_{f_i}L_{f_i}}{N^2}\right)
					\|x_i-y_i\|^2
					{}-{}
					\tfrac{\gamma_i\mu_{f_i}}{N}
					\left(2-\tfrac{\gamma_i}{N}(\mu_{f_i}+L_{f_i})\right)
					\|x_i-y_i\|^2
				\\
				{}={} &
					\left(1-\tfrac{\gamma_i\mu_{f_i}}{N}\right)^2
					\|x_i-y_i\|^2,
				\end{align*}
				where strong convexity and the fact that $\gamma_i<\nicefrac{N}{L_{f_i}}\leq\nicefrac{2N}{(\mu_{f_i}+L_{f_i})}$ was used in the second inequality.
				Multiplying by $\gamma_i^{-1}$ and summing over $i$ shows that \(\id-\Gamma\nabla F\) is \((1-\delta)\)-contractive in the metric \(\|{}\cdot{}\|_{\Gamma^{-1}}\), and so is \(\T=\prox_G^{\Gamma^{-1}}\circ(\Fw{})\) as it follows from assert \ref{thm:FNE}.%
			\qedhere
			\end{proofitemize}
		\end{proof}
	\end{lem}
	The next result recaps an important property that the FBE inherits from the cost function \(\Phi\) that is instrumental for establishing global convergence and asymptotic linear rates for the BC-\Cref{alg:BC}.
	The result falls as special case of \cite[Th. 5.2]{yu2019deducing} after observing that
	\[
		\FBE(\bm x)
	{}={}
		\inf_{\bm w}\set{
			\Phi(\bm w)
			{}+{}
			D_H(\bm w,\bm x)
		},
	\]
	where
	\(
		D_H(\bm w,\bm x)
	{}={}
		H(\bm w)-H(\bm x)-\innprod{\nabla H(\bm x)}{\bm w-\bm x}
	\)
	is the Bregman distance with kernel \(H=\tfrac12\|{}\cdot{}\|_{\Gamma^{-1}}^2-F\).
	
	\begin{lem}[{\cite[Th. 5.2]{yu2019deducing}}]\label{thm:loja}%
		Suppose that \Cref{ass:basic} holds and for \(\gamma_i\in(0,\nicefrac{N}{L_{f_i}})\), \(i\in[N]\), let \(\Gamma=\blockdiag(\gamma_1\I_{n_1},\dots,\gamma_N\I_{n_N})\).
		If \(\Phi\) has the KL property with exponent \(\theta\in(0,1)\) (as is the case when \(f_i\) and \(G\) are semialgebraic), then so does \(\FBE\) with exponent
		\(
			\max\set{\nicefrac12,\theta}
		\).
	\end{lem}
	
	\ifaccel
		\begin{lem}[FBE: convexity and block-smoothness]\label{thm:convex}%
			Suppose that \Cref{ass:basic,ass:Fast} are satisfied, and consider the notation introduced therein.
			Let \(\gamma_i\in(0,\nicefrac{N}{L_{f_i}})\) be fixed.
			Define
			\(
				Q_i
			{}\coloneqq{}
				\gamma_i^{-1}\I-\tfrac{1}{N}H_i\in\R^{n_i\times n_i}
			\),
			\(
				Q
			{}\coloneqq{}
				\blockdiag(Q_1,\dots,Q_N)
			\),
			and
			\(
				H
			{}\coloneqq{}
				\tfrac1N\blockdiag(H_1,\dots,H_N)
			\).
			Then, $\FBEC \coloneqq\FBE \circ Q^{-1/2}$ is convex and smooth with $\nabla\FBEC (\tilde{\bm x}) = Q^{1/2}(\bm x-\T(\bm x))$ where $\bm x=Q^{-1/2}\tilde{\bm x}$. 
			In fact, for any \(\tilde{\bm x},\tilde{\bm x}'\in\R^{\sum_in_i}\) it holds that
			\begin{equation}\label{eq:FBEComposedsmooth}
				0
			{}\leq{}
				\innprod{\nabla\FBEC(\tilde{\bm x}')-\nabla\FBEC(\tilde{\bm x})}{\tilde{\bm x}'-\tilde{\bm x}}
			{}\leq{}
				\|\tilde{\bm x}'-\tilde{\bm x}\|^2.
			\end{equation}
			In particular, function $\FBEC$ is $1$-smooth along each block $i\in[N]$.  	If, additionally, all functions \(f_i\) are strongly convex, then 
			$\FBEC$ is \(\sigma\)-strongly convex with $\sigma\coloneqq \tfrac{1}{N}\min_{i\in [N]}\left\{\gamma_i\mu_{f_i}\right\}$. 
			\begin{proof}
				Since $\gamma_i<N/L_{f_i}$, $Q$ is positive definite.
				We begin by showing that for any \(\bm x,\bm x'\in\R^{\sum_in_i}\) it holds that
				\begin{equation}\label{eq:FBEsmooth}
					0\leq \|\bm x'-\bm x\|^2_{Q}- \|Q(\bm x'-\bm x)\|^2_{\Gamma}
				{}\leq{}
					\innprod{\nabla\FBE(\bm x')-\nabla\FBE(\bm x)}{\bm x'-\bm x}
				{}\leq{}
					\|\bm x'-\bm x\|^2_Q.
				\end{equation}
				It follows from \Cref{thm:MoreauGrad}, the chain rule of differentiation applied to \eqref{eq:FBEMoreau}, and the twice continuous differentiability of \(F\) that \(\FBE\) is continuously differentiable with
				\(
					\nabla\FBE(\bm x)
				{}={}
					Q(\bm x-\bm z)
				\).
				For \(\bm z^x\coloneqq\T(\bm x)\) and \(\bm z^{x'}\coloneqq\T(\bm {x'})\) it holds that
				\begin{equation}\label{eq:innprodGrad}
					\innprod{\nabla\FBE(\bm {x'})-\nabla\FBE(\bm x)}{\bm {x'}-\bm x}
				{}={}
					\innprod{
						Q(\bm {x'}-\bm z^{x'}-\bm x+\bm z^x)
					}{
						\bm {x'}-\bm x
					}
				{}={}
					\|\bm {x'}-\bm x\|^2_Q
					{}-{}
					\innprod{
						\bm z^{x'}-\bm z^x
					}{
						Q(\bm {x'}-\bm x)
					}.
				\end{equation}
				In order to bound the last scalar product, observe that
				\[
					0
				{}\leq{}
					\innprod{
						\Gamma^{-1}(\bm z^{x'}-\bm z^x)
					}{
						(\bm {x'}-\Gamma\nabla F(\bm {x'}))
						{}-{}
						(\bm x-\Gamma\nabla F(\bm x))
					}
				{}\leq{}
					\bigl\|
						(\bm {x'}-\Gamma\nabla F(\bm {x'}))
						{}-{}
						(\bm x-\Gamma\nabla F(\bm x))
					\bigr\|_{\Gamma^{-1}}^2,
				\]
				as it follows from \Cref{thm:FNE}.
				Since \(\id-\Gamma\nabla F=\Gamma Q{}\cdot{} - \Gamma\bm q\) (with $\bm q\coloneqq(\tfrac1Nq_1,\dots,\tfrac1Nq_N)$), the above inequality simplifies to
				\[
					0
				{}\leq{}
					\innprod{
						\bm z^{x'}-\bm z^x
					}{
						Q(\bm {x'}-\bm x)
					}
				{}\leq{}
					\|\Gamma Q(\bm {x'}-\bm x)\|_{\Gamma^{-1}}^2,
				\]
				which combined with \eqref{eq:innprodGrad} results in the claimed \eqref{eq:FBEsmooth}.
				If additionally \(\mu_{f_i}>0\) for all \(i\), then \(\FBE\) is \(1\)-strongly convex in the metric \(\|{}\cdot{}\|^2_{Q - Q\Gamma Q}\) (by observing that \(Q-Q\Gamma Q\succ 0\)). The result in \eqref{eq:FBEComposedsmooth} follows by using \eqref{eq:FBEsmooth} with the change of variables $\bm x=Q^{-1/2}\tilde{\bm x}$, $\bm x'=Q^{-1/2}\tilde{\bm x}'$ and noting that $\nabla \FBEC(\tilde{\bm x}) = Q^{-1/2}\nabla \FBE (\bm x)$. 
				Since $\Gamma$ is block-wise a multiple of identity it commutes with any block-diagonal matrix. Therefore, when $f_i$ are strongly convex, using the lower bound in \eqref{eq:FBEsmooth} and the above change of variable we obtain that  $\FBEC$ is strongly convex in the metric \(\|{}\cdot{}\|^2_{\I - \Gamma Q}\). The result follows by noting that $\I - \Gamma Q= \Gamma H$. 
			\end{proof}
		\end{lem}
	\fi

	\end{appendix}


	\ifarxiv
		\bibliographystyle{plain}
	\else
		\phantomsection
		\addcontentsline{toc}{section}{References}
		\bibliographystyle{spmpsci}
	\fi
	\bibliography{Bibliography.bib}

\end{document}